\documentclass[12pt]{article}
\usepackage{srcltx}
\usepackage{amssymb}
\usepackage{amsfonts}
\usepackage{amsmath}
\usepackage{amsthm}
\usepackage{enumerate}
\usepackage{multicol}

\newtheorem{theorem}{Theorem}[section]
\newtheorem{prop}[theorem]{ Proposition}
\newtheorem{lemma}[theorem]{Lemma}

\newtheorem{remark}{Remark}[section]

\newcommand\cA{{\cal A}}
\newcommand\cC{{\cal C}}
\newcommand\cE{{\cal E}}

\newcommand\cL{{\cal L}}
\newcommand\cB{{\cal B}}

\newcommand\cN{{\cal N}}

\newcommand\cP{{\cal P}}

\newcommand\cR{{\cal R}}
\newcommand\cT{{\cal T}}

\newcommand\e{\epsilon}
\newcommand\ve{\varepsilon}

\newcommand\ov{\overline}

\def\bbr{{\mathbb R}}

\def\text#1{\hbox{#1}}

\def\endproof{\mbox{\ $\qed$}}

\def\E{{\bf E}}
\def\e{{\bf e}}
\def\A{{\bf A}}
\def\P{{\bf P}}

\def\H{{\bf H}}
\def\L{{\bf L}}

\newcommand{\wh}{\widehat}
\newcommand{\wt}{\widetilde}

\def\O{\hbox{\rm O}}

\def\Chi{{\bf 1}}
\def\d{\mathrm{d}}
\def\build #1_#2{\mathrel{\mathop{\kern 0pt #1}\limits_{#2}}} 
\newcommand{\zs}[1]{{\mathchoice{#1}{#1}{\lower.25ex\hbox{$\scriptstyle#1$}}
{\lower0.25ex\hbox{$\scriptscriptstyle#1$}}}}
\numberwithin{equation}{section}

\begin{document}
\title{Adaptive asymptotically efficient estimation\\
%{\Large }  
in heteroscedastic nonparametric regression.
}
\author{{\Large By  Leonid Galtchouk and Sergey Pergamenshchikov}\\
%\thanks{The second author is partially supported by the RFFI-Grant 04-01-00855.}\\
Louis Pasteur University  of Strasbourg and University of Rouen
}
\date{}
\maketitle

\begin{abstract}
The paper deals with asymptotic properties of the adaptive
procedure  proposed in the author paper, 2007, 
for  estimating an unknown nonparametric regression.
%\cite{GaPe1}. 
We prove that this
procedure is asymptotically efficient for a quadratic risk,
i.e. the asymptotic quadratic risk for this procedure
coincides with the Pinsker constant which gives a sharp 
 lower bound for the quadratic risk over all possible estimators.
\footnote{
{\sl AMS 2000 Subject Classification} : primary 62G08; secondary 62G05, 62G20}
\footnote{
{\sl Key words}: asymptotic bounds, adaptive estimation, efficient estimation,
heteroscedastic regression, nonparametric regression, Pinsker's constant.
}
\end{abstract}
\bibliographystyle{plain}
\renewcommand{\columnseprule}{.1pt}
\newpage

\section{Introduction}\label{sec:I}

The paper deals with the  estimation problem in 
the heteroscedastic nonparametic regression model 
\begin{equation}\label{sec:I.1}
y_\zs{j}=S(x_\zs{j})+\sigma_\zs{j}(S)\,\xi_\zs{j}\,, 
\end{equation}
where the design points $x_\zs{j}=j/n$, $S(\cdot)$
 is an unknown function to be estimated, 
$(\xi_\zs{j})_\zs{1\le j\le n}$ is a sequence of centered 
%i.i.d.
 independent random variables with
unit variance and 
$(\sigma_\zs{j}(S))_\zs{1\le j\le n}$ are unknown scale functionals depending on
the design points and the regression function $S$.

 Typically, the notion of asymptotic
optimality is associated with the optimal convergence rate of the minimax risk
(see e.g., Ibragimov, Hasminskii,1981; Stone,1982).
% \cite{sec:IbHa}, \cite{St})
  An important question in optimality
results is to study the exact asymptotic behavior of the minimax risk. Such results
have been obtained only in a limited number of investigations. As to the
nonparametric estimation problem
for heteroscedastic regression models we should mention the papers by Efroimovich, 2007,
Efroimovich, Pinsker, 1996, and Galtchouk, Pergamenshchikov, 2005,
% \cite{Ef}-\cite{EfPi} and \cite{GaPe0} 
 concerning the exact asymptotic behavior of the $\cL_\zs{2}$-risk
and the paper by Brua, 2007,
%\cite{Br}  
devoted to the efficient pointwise estimation for
heteroscedastic regressions. 

Heteroscedastic regression models are extensively used in the financial mathematics, in particular,
 in problems of calibrating  (see e.g., Belomestny, Reiss, 2006).
Moreover, these models are  popular in econometrics
(see, for e.g., Goldfeld, Quandt, 1972,
%\cite{GoQu},
 p. 83), which, for exemple,
for  consumer budget problems, makes use of
  some semiparametric version of 
model \eqref{sec:I.1}  with the scale coefficients  of type
\begin{equation}\label{sec:I.2}
\sigma^2_\zs{j}(S)=c_\zs{0}+c_\zs{1}x_\zs{j}+c_\zs{2}S^2(x_\zs{j})+
c_\zs{3} \int^1_\zs{0}\,S^2(t) \d t
\,,
\end{equation}
where $(c_\zs{i})_\zs{0\le i\le 3}$ are some unknown positive  constants.

The goal of this paper is to study asymptotic properties of the adaptive 
estimation procedure proposed in Galtchouk, Pergamenshchikov, 2007,
%\cite{GaPe1} 
for which a non-asymptotic oracle
inequality was proved for quadratic risks. More precisely,
in this paper
we show that this procedure is efficient under some conditions on the scale functions
$(\sigma_\zs{j}(S))_\zs{1\le j\le n}$ which hold for the functions \eqref{sec:I.2}.
Note that in  Efroimovich, 2007, Efroimovich, Pinsker, 1996,
%\cite{Ef}-\cite {EfPi} 
an efficient adaptive procedure is constructed
for heteroscedastic regression when the scale coefficient is independent of $S$, i.e.
$\sigma_\zs{j}(S)=\sigma_\zs{j}$. In Galtchouk, Pergamenshchikov, 2005,
%\cite{GaPe0} 
for the model \eqref{sec:I.1}, the asymptotic efficiency
was proved under conditions
 which are not satisfied in the case \eqref{sec:I.2}.
Moreover, in the these papers the efficiency is proved only for the gaussian random variables
$(\xi_\zs{j})_\zs{1\le j\le n}$ in model \eqref{sec:I.1}. This is  very restrictive condition
 for applications.

In this paper we consider the robust quadratic risk, i.e.
in the definition of the risk
 we take the additional supremum over the family
of unknown noise distributions likely to Galtchouk, Pergamenshchikov, 2006.
%\cite{GaPe}
 This modification allows us to eliminate from
the risk dependence on the noise distribution. Moreover, for this risk the efficient procedure 
is robust with respect to changing the noise distribution. 

As is well known, to prove the asymptotic efficiency one has to show that the asymptotic quadratic 
risk coincides with the lower bound which is equal to the Pinsker constant. In the paper two 
problems
are resolved: in the first one a upper bound for the risk is obtained by making use of the
non-asymptotic oracle inequality from Galtchouk, Pergamenshchikov, 2007,
%\cite{GaPe1}
in the second one we prove that this upper
bound coincides with the Pinsker constant. Let us remember that the adaptive procedure proposed in
%\cite{GaPe1} 
Galtchouk, Pergamenshchikov, 2007, is based on weighted least squares estimators, where the 
weights are proper
modifications of the Pinsker weights for the homogeneous case (when
$\sigma_1(S)=\ldots=\sigma_\zs{n}(S)=1$) relative to a certain smoothness of the function $S$ and
this procedure chooses a best estimator for the quadratic risk among these estimators. To obtain 
the Pinsker constant for the model \eqref{sec:I.1} one has to prove a sharp asymptotic lower bound 
for the quadratic risk in the case when the noise variance depends on the unknown regression
function.
In this case, as usually, we minorize the minimax risk by a bayesian one for a respective parametric family. Then for the bayesian risk we make use of a lower bound (see Theorem 6.1)
which is a modification of the van Trees inequality (see, Gill, Levit, 1995).

 The paper is organized as follows.  
In Section~\ref{sec:Ad} we construct an adaptive estimation procedure.
Section~\ref{sec:Co} contains the principal conditions. The statements of our major results (the upper and lower bounds for quadratic risks)
are presented in Section~\ref{sec:Ma}. The upper bound is proved
in Section~\ref{sec:Up}. 
In Section~\ref{sec:Lo} we give all main steps of proving the lower bound:
in Subsection~\ref{subsec:Tr} we find the lower bound for the bayesian
 risk for a parametric regression model which minorizes the minimax risk;
in Subsection~\ref{subsec:Fa} we study a special family of parametric functions  
used to define the bayesian risk;
in Subsection~\ref{subsec:Br} we choose a prior distribution for bayesian risk
to maximize the lower bound.
Section~\ref{sec:Np} is devoted to explain how to use the given procedure
 in the case, when the unknown regression function is non periodic.
 In Section~\ref{sec:Cn} we discuss the main results and their practical importance.
The proofs are given in Section~\ref{sec:Pr}.

\section{Adaptive procedure}\label{sec:Ad}

In this section we describe the adaptive procedure proposed in 
Galtchouk, Pergamenshchikov, 2006. 
To evaluate the error of estimation in the model \eqref{sec:I.1} we make use of the
empiric quadratic norm in the Hilbert space $\cL_\zs{2}[0,1]$
generated by the design points $(x_\zs{j})_\zs{1\le j\le n}$ of model \eqref{sec:I.1}. To this
end, for any functions $u$ and $v$ from $\cL_\zs{2}[0,1]$, we define the empiric
inner product
\begin{equation}\label{sec:Ad.0}
(u,v)_\zs{n}=
\frac{1}{n}\sum^n_\zs{l=1}\,u(x_\zs{l})\,v(x_\zs{l})\,.
\end{equation}
Therefore, the  estimation error of an estimator $\wh{S}$ of $S$ will be evaluated
by the empiric quadratic loss function
$$
\|\wh{S}-S\|^2_\zs{n}=\frac{1}{n}\sum^n_\zs{l=1}\,(\wh{S}(x_\zs{l})-S(x_\zs{l}))^2\,.
$$
Moreover, we make use of  this inner product for vectors in $\bbr^n$ as well, i.e.
if
$u=(u_\zs{1},\ldots,u_\zs{n})'$ and  $v=(v_\zs{1},\ldots,v_\zs{n})'$, then 
$$
(u,v)_\zs{n}=\frac{1}{n}u'v=
\frac{1}{n}\sum^n_\zs{l=1}\,u_\zs{l}\,v_\zs{l}\,.
$$
The prime denotes the transposition.

 Let now 
$(\phi_\zs{j})_\zs{j\ge 1}$
be 
 the 
 standard trigonometric basis in $\cL_2[0,1]$, i.e.
\begin{equation}\label{sec:Ad.0-0}
\phi_1(x)=1\,,\quad
\phi_\zs{j}(x)=\sqrt{2}\,Tr_\zs{j}(2\pi [j/2]x)\,,\ j\ge 2\,,
\end{equation}
where the function $Tr_\zs{j}(x)=\cos(x)$ for even $j$ and
$Tr_\zs{j}(x)=\sin(x)$ for odd $j$; $[x]$ denotes the integer part of $x$.

Notice that if $n$ is odd, then the functions $(\phi_\zs{j})_\zs{1\le j\le n}$ 
are orthonormal with respect to  the empiric inner product \eqref{sec:Ad.0},
 i.e. for any $1\le i,j\le n$,
\begin{equation}\label{sec:Ad.0-1}
(\phi_\zs{i}\,,\,\phi_\zs{j})_\zs{n}=
\frac{1}{n}\sum^n_\zs{l=1}\phi_\zs{i}(x_l)\phi_\zs{j}(x_l)={\bf Kr}_\zs{ij}\,,
\end{equation}
where ${\bf Kr}_\zs{ij}$ is Kronecker's symbol,  ${\bf Kr}_\zs{ij}=1$ if $i=j$ and 
 ${\bf Kr}_\zs{ij}=0$ for $i\ne j$.

\begin{remark}\label{Re.Ad.1}
Note that in the case of even $n$, the basis  \eqref{sec:Ad.0-0} is orthogonal and it is orthonormal except
of the  $n$th function for which the normalizing constant should be changed.
The corresponding modifications, for even $n$,  one can see in Galtchouk, Pergamenshchikov,2005.
To avoid these complications, we suppose $n$ to be odd.
\end{remark}

Thanks to this basis we pass to  the discrete Fourier transformation of model \eqref{sec:I.1}:
\begin{equation}\label{sec:Ad.1}
\wh{\theta}_\zs{j,n}=\theta_\zs{j,n}+\frac{1}{\sqrt{n}}\xi_\zs{j,n}\,,
\end{equation}
where $\wh{\theta}_\zs{j,n}=(Y,\phi_\zs{j})_\zs{n}$, $Y=(y_\zs{1},\ldots,y_\zs{n})'$,
$\theta_\zs{j,n}=(S,\phi_\zs{j})_\zs{n}$
and
$$
\xi_\zs{j,n}=\frac{1}{\sqrt{n}}\sum^n_\zs{l=1}\sigma_l(S)\xi_l\phi_\zs{j}(x_l)\,.
$$
We estimate the function $S$ by the weighted least squares estimator
\begin{equation}\label{sec:Ad.2}
\wh{S}_\zs{\lambda}=\sum^n_\zs{j=1}\lambda(j)\wh{\theta}_\zs{j,n}\phi_\zs{j}\,,
\end{equation}
where the weight vector $\lambda=(\lambda(1),\ldots,\lambda(n))'$
belongs to some finite set $\Lambda$ from $[0,1]^n$ with $n\ge 3$.

Here we make use of the weight family $\Lambda$ introduced
in Galtchouk, Pergamenshchikov, 2009, i.e.
\begin{equation}\label{sec:Ad.3}
\Lambda\,=\,\{\lambda_\zs{\alpha}\,,\,\alpha\in\cA\}\,,
\quad
\cA=\{1,\ldots,k^*\}\times\{t_1,\ldots,t_m\}\,,
\end{equation}
where  $t_\zs{i}=i\ve$ and $m=[1/\ve^2]$. We suppose that the parameters 
$k^*\ge 1$ and $0<\ve\le 1$ are functions of $n$, i.e. $k^*= k^*_\zs{n}$ and 
$\ve=\ve_\zs{n}$, such that,
\begin{equation}\label{sec:Ad.3-1}
\left\{
\begin{array}{ll}
&\lim_\zs{n\to\infty}\,k^*_\zs{n}=+\infty\,,
\quad
\lim_\zs{n\to\infty}\,\dfrac{k^*_\zs{n}}{\ln n}=0\,,\\[4mm]
&\lim_\zs{n\to\infty}\,\ve_\zs{n}\,=\,0
\quad\mbox{and}\quad
\lim_\zs{n\to\infty}\,n^{\nu}\,\ve_\zs{n}\,=+\infty\,,
\end{array}
\right.
\end{equation}
for any $\nu>0$. For example, one can take for  $n\ge 3$
$$
\ve_\zs{n}=1/\ln n
\quad\mbox{and}\quad
k^*_\zs{n}=\ov{k}+\sqrt{\ln n}\,,
$$
where $\ov{k}$ is any nonnegative constant.

For each $\alpha=(\beta,t)\in\cA$ we define the weight vector
$\lambda_\zs{\alpha}=(\lambda_\zs{\alpha}(1),\ldots,\lambda_\zs{\alpha}(n))'$ as
\begin{equation}\label{sec:Ad.4}
\lambda_\zs{\alpha}(j)=\Chi_\zs{\{1\le j\le j_\zs{0}\}}+
\left(1-(j/\omega(\alpha))^\beta\right)\,
\Chi_\zs{\{ j_\zs{0}<j\le \omega(\alpha)\}}\,.
\end{equation}
Here $j_0=j_\zs{0}(\alpha)=\left[\omega(\alpha)\,\ve_\zs{n}\right]$ with
\begin{equation}\label{sec:Ad.4-1}
\omega(\alpha)=\ov{\omega}+(A_\zs{\beta}\,t)^{1/(2\beta+1)}n^{1/(2\beta+1)}\,,
 \end{equation}
where $\ov{\omega}$ is any nonnegative constant 
and
$$
A_\zs{\beta}=\frac{(\beta+1)(2\beta+1)}{\beta\pi^{2\beta}}\,.
$$

\begin{remark}\label{Re.Ad.2}
Note that,  the weighted least squares estimators 
 \eqref{sec:Ad.2} have been introduced by Pinsker, 1981, for the optimal filtering of a continuous time  signal  in the gaussian noise. It tourned out, that the asymptotic quadratic risk 
for  estimators of type \eqref{sec:Ad.2} -- \eqref{sec:Ad.4} is minimal over all possible estimators.
This sharp  value of the asymptotic
 quadratic risk is called  the Pinsker constant.
 Nussbaum, 1985, makes
 use of the same method with proper modification for efficient
estimation of the function $S$ of known smoothness
in the homogeneous gaussian model \eqref{sec:I.1}, i.e. when
$\sigma_1(S)=\ldots=\sigma_\zs{n}(S)=1$ and $(\xi_\zs{j})_\zs{1\le j\le n}$ is a sequence of i.i.d. $\cN(0,1)$ random variables. 
\end{remark}
 To choose weights from the set \eqref{sec:Ad.3} we minimize
the special cost function introduced by
 Galtchouk, Pergamenshchikov, 2007. This cost function is as follows
\begin{equation}\label{sec:Ad.5}
J_\zs{n}(\lambda)\,=\,\sum^n_\zs{j=1}\,\lambda^2(j)\wh{\theta}^2_\zs{j,n}\,-
2\,\sum^n_\zs{j=1}\,\lambda(j)\,\wt{\theta}_\zs{j,n}\,
+\,\rho \wh{P}_\zs{n}(\lambda)\,,
\end{equation}
where 
\begin{equation}\label{sec:Ad.6}
\wt{\theta}_\zs{j,n}=
\wh{\theta}^2_\zs{j,n}-\frac{1}{n}\wh{\varsigma}_\zs{n}
\quad\mbox{with}\quad
\wh{\varsigma}_\zs{n}=\sum^{n}_\zs{j=l_\zs{n}+1}
\wh{\theta}^2_\zs{j,n}
\end{equation}
and $l_\zs{n}=[n^{1/3}+1]$. The penalty term we define as
$$
\wh{P}_\zs{n}(\lambda)=\frac{|\lambda|^2 \wh{\varsigma}_\zs{n}}{n}\,,\quad
|\lambda|^2=\sum^n_\zs{j=1} \lambda^2(j)
\quad\mbox{and}\quad
\rho=\frac{1}{3+L_\zs{n}}\,,
$$
where $L_\zs{n}\ge 0$ is any slowly increasing sequence, i.e.
\begin{equation}\label{sec:Ad.6-1}
\lim_\zs{n\to\infty}\,L_\zs{n}=+\infty
\quad\mbox{and}\quad
\lim_\zs{n\to\infty}\,\frac{L_\zs{n}}{n^\nu}=0\,,
\end{equation}
for any $\nu>0$.

Finally, we set
\begin{equation}\label{sec:Ad.7}
\wh{\lambda}=\mbox{argmin}_\zs{\lambda\in\Lambda}\,J_\zs{n}(\lambda)
\quad\mbox{and}\quad
\wh{S}_\zs{*}=\wh{S}_\zs{\wh{\lambda}}\,.
\end{equation}

The goal of this paper is to study asymptotic properties (as $n\to\infty$)
of this estimation procedure.

\begin{remark}\label{Re.Ad.3}
Now we explain why does one choose the cost function in the form  \eqref{sec:Ad.5}.
Developing the empiric quadratic loss function for the estimator \eqref{sec:Ad.2}, one obtains
$$
\|\wh{S}_\zs{\lambda}-S\|^2_\zs{n}=\sum^n_\zs{j=1}\,\lambda^2(j)\wh{\theta}^2_\zs{j,n}\,-
2\,\sum^n_\zs{j=1}\,\lambda(j)\,\wh{\theta}_\zs{j,n}\,\theta_\zs{j,n} +\|S\|^2_\zs{n}\,.
$$
It's natural to choose the weight vector  $\lambda$ for which this function reaches
the minimum. Since the last term on the right-hand part is independent of  $\lambda$, it can be dropped and one has to minimize with respect to $\lambda$ the function equals to the difference of the two first terms on the right-hand part. It's clear that the minimization problem can not be solved directly because the Fourier coefficients $(\theta_\zs{j,n}) $ are unknown. To overcome
this difficulty, we replace the product $\wh{\theta}_\zs{j,n}\,\theta_\zs{j,n}$
by its asymptotically unbiased estimator 
$\wt{\theta}_\zs{j,n}$. Moreover, to pay this
substitution, we introduce into the cost function the penalty term $\wh{P}_\zs{n}$ with a small coefficient $\rho>0$. The form of the penalty term is provided by the principal term of the quadratic risk for weighted least-squares estimator,
see Galtchouk, Pergamenshchikov, 2007, 2009. The coefficient $\rho>0$ is small, because the estimator
 $\wt{\theta}_\zs{j,n}$ is close in mean to the quantity
$\wh{\theta}_\zs{j,n}\,\theta_\zs{j,n}$ asymptotically, as $n\to\infty$.

Note that the principal difference between the procedure \eqref{sec:Ad.7} and the adaptive procedure proposed by Golubev, Nussbaum, 1993, for a homogeneous gaussian regression, consists in 
presence of the penalty term in the cost function \eqref{sec:Ad.5}. 
\end{remark}

\begin{remark}\label{Re.Ad.4}
As it was noted in Remark~\ref{Re.Ad.2}, Nussbaum, 1985, has shown that the weight coefficients
of type \eqref{sec:Ad.4} provide the asymptotic minimum of the quadratic risk in the
estimation problem of the regression function  for the homogeneous gaussian model  \eqref{sec:I.1}, when the smoothness of the function $S$ is known. In fact, to obtain an efficient estimator one needs to take a weighted least squares estimator \eqref{sec:Ad.2} with the weight vector
 $\lambda_\zs{\alpha}$, where the index $\alpha$ depends on the smoothness of the function $S$
(the parameters $k$ and $r$ in \eqref{sec:Co.1})
 and on the
 coefficients $(\sigma_\zs{j}(S))_\zs{1\le j\le n}$,
(the paramter  $\varsigma(S)$ in \eqref{sec:Co.5}), which are unknown in our case. For this reason, Galtchouk,
 Pergamenshchikov , 2007, have proposed to make use of the family of coefficients \eqref{sec:Ad.3}, which contains 
the   weight vector providing the minimum of the quadratic
risk. Indeed, the set
  \eqref{sec:Ad.3} is a two-dimensional grille giving all possible
values for
 the  parameters $k$ and 
$\ov{r}(S)=r/\varsigma(S)$ in the  estimator \eqref{sec:Ad.2} with
the weight
\eqref{sec:Ad.4} which is efficient if the parameters $k$ and $\ov{r}(S)$ are known, i.e.
the set $\cA$ gives the familly of efficient estimators.
Moreover, under some weak conditions on the coefficients
$(\sigma_\zs{j}(S))_\zs{1\le j\le n}$,
Galtchouk,
 Pergamenshchikov , 2007, shown that
the
procedure  \eqref{sec:Ad.7} is best in the class of these estimators
in the sens of the 
non-asymptotic oracle inequality (see, Theorem~\ref{Th.sec:Ma.1} below).
It is important to note that, due to the conditions
\eqref{sec:Ad.3-1} for the set $\cA$, the secondary term in the oracle inequality is slowly increasing (slower than any degree of $n$).
\end{remark}

\section{Conditions}\label{sec:Co}

First we impose some conditions on the function $S$ in the model \eqref{sec:I.1}.

Let   $\cC^{k}_\zs{per,1}(\bbr)$ be  the set of $1$-periodic
$k$ times differentiable $\bbr\to\bbr $ functions. We assume that
 $S$ belongs to the following set:
\begin{equation}\label{sec:Co.1}
W^{k}_\zs{r}=\{f\in\cC^{k}_\zs{per,1}(\bbr)
\,:\,\sum_\zs{j=0}^k\,\|f^{(j)}\|^2\le r\}\,,
 \end{equation}
where $\|\cdot\|$ denotes the  norm in $\cL_\zs{2}[0,1]$, i.e.
\begin{equation}\label{sec:Co.2}
\|f\|^2=\int^1_\zs{0}f^2(t)\d t\,.
\end{equation}
Moreover, we suppose that $r>0$ and $k\ge 1$ are unknown parameters.

Note that the set $W^{k}_\zs{r}$ can be represented as
an ellipse in $\cL_\zs{2}[0,1]$, i.e.
\begin{equation}\label{sec:Co.3}
W^{k}_\zs{r}=\{f\in\cL_\zs{2}[0,1]\,:\,
\sum_\zs{j=1}^\infty\,a_\zs{j}\theta^2_\zs{j}\le r\}\,,
 \end{equation}
where 
\begin{equation}\label{sec:Co.3-1}
\theta_\zs{j}=(f,\phi_\zs{j})=\int^1_\zs{0}f(t)\phi_\zs{j}(t)\d t
\end{equation}
and 
\begin{equation}\label{sec:Co.4}
a_\zs{j}=\sum^k_\zs{l=0}\|\phi^{(l)}_\zs{j}\|^2=
\sum^k_{i=0}(2\pi [j/2])^{2i}\,.
\end{equation}
Here $(\phi_\zs{j})_\zs{j\ge 1}$ is the trigonometric basis 
defined in \eqref{sec:Ad.0-0}.

Now we describe the conditions on the scale coefficients $(\sigma_\zs{j}(S))_\zs{j\ge 1}$.

\begin{itemize}
\item[$\H_\zs{1})$] {\em 
 $\sigma_\zs{j}(S)=g(x_\zs{j},S)$ for some unknown  function 
$g : [0,1]\times \cL_\zs{1}[0,1] \to \bbr_+$, which is
square integrable with respect to $x$ and such that
\begin{equation}\label{sec:Co.5}
\lim_\zs{n\to\infty}\,\sup_\zs{S\in W^k_r}\,\left|
\,\frac{1}{n}\,\sum^n_\zs{j=1}\,g^2(x_\zs{j},S)\,-\,\varsigma(S)\,\right|\,=0\,,
\end{equation}
where $\varsigma(S):=\,\int_0^1\,g^2(x,S)\d x$. 
Moreover,
\begin{equation}\label{sec:Co.6}
g_\zs{*}=\inf_\zs{0\le x\le 1}\,
\inf_\zs{S\in W^k_r} g^2(x,S)\,>0
\end{equation}
and
\begin{equation}\label{sec:Co.6-1}
\sup_\zs{S\in W^k_r} \varsigma(S)<\infty\,.
\end{equation}
}
\item[$\H_\zs{2})$] {\em 
For any $x\in [0,1]$, the operator 
$g^2(x,\cdot)\,:\,\cC[0,1]\to \bbr$
 is differentiable in  
the Fr\'echet sense at any fixed function $f_\zs{0}$ from $\cC[0,1]$
, i.e.
 for any $f$ from some vicinity of $f_\zs{0}$ in $\cC[0,1]$,
$$
g^2(x,f)=g^2(x,f_\zs{0})+\L_\zs{x,f_\zs{0}}(f-f_\zs{0})+
\Upsilon(x,f_\zs{0},f)\,,
$$
where the Fr\'echet derivative
$\L_\zs{x,f_\zs{0}}\,:\,\cC[0,1]\to \bbr$
is a bounded linear operator
and 
the residual term $\Upsilon(x,f_\zs{0},f)$, for each $x\in [0,1]$, satisfies the following
property:
$$
\lim_\zs{\|f-f_\zs{0}\|_\zs{\infty}\to 0}
\frac{|\Upsilon(x,f_\zs{0},f)|}{\|f-f_\zs{0}\|_\zs{\infty}}=0\,,
$$
where $\|f\|_\zs{\infty}=\sup_\zs{0\le t\le 1} |f(t)|$.
}

\item[$\H_\zs{3})$] {\em
There exists some positive constant $C^*$ such that,
for any function $S$ from $\cC[0,1]$, the operator 
$\L_\zs{x,S}$ defined in the condition $\H_\zs{2})$ 
satisfies the following inequality, for any function $f$ from $\cC[0,1]$:
\begin{equation}\label{sec:Co.8}
|\L_\zs{x,S}(f)|
\le C^*
\left(
|S(x)f(x)|+|f|_\zs{1}+\|S\|\,\|f\|
\right)\,,
\end{equation}
where $|f|_\zs{1}=\int^1_\zs{0}|f(t)|\d t$.
}

\item[$\H_\zs{4})$] {\em The function 
$g_\zs{0}(\cdot)=g(\cdot,S_\zs{0})$ corresponding to $S_\zs{0}\equiv 0$
is  continuous on the interval $[0,1]$.
Moreover,
$$
\lim_\zs{\delta\to 0}\,
\sup_\zs{0\le x\le 1}\,
\sup_\zs{\|S\|_\zs{\infty}\le \delta}\,
|g(x,S)-g(x,S_\zs{0})|\,=\,0\,.
$$
}
\end{itemize}

\begin{remark}\label{Re.sec:Co.1}
Let us explain the conditions $\H_\zs{1})$-$\H_\zs{4})$ which are the regularity conditions of the function $g(x,S)$ generating the scale coefficients $(\sigma_\zs{j}(S))_\zs{1\le j\le n}$.

Condition $\H_\zs{1})$ means that the function $g(\cdot,S)$ must be uniformly integrable
with respect to the first argument in the sense of convergence \eqref{sec:Co.5}.
Moreover, this function must be separated  from zero (see inequality \eqref{sec:Co.6})
and bounded on the class \eqref{sec:Co.1} (see inequality \eqref{sec:Co.6-1}). Boundedness away
from zero provides that the distribution of observations $(y_\zs{j})_\zs{1\le j\le n}$ isn't
degenerate in $\bbr^n$, and the boundedness means that the intensity of the noise vector must be finite, otherwise the estimation problem has not any sense.
 
 Conditions $\H_\zs{2})$ and $\H_\zs{3})$ mean that the function
$g(x,\cdot)$ is regular, at any fixed $0\le x\le 1$,  with respect to $S$ in the sense, that it is
differentiable in the Fr\`echet sense (see e.g., Kolmogorov, Fomin, 1989) and, moreover, the Fr\`echet derivative satisfies the growth condition given by the inequality \eqref{sec:Co.8}
 which permits to consider the example \eqref{sec:I.2}.
 
The last condition $\H_\zs{4})$ is the usual uniform continuity condition of the function $g(\cdot,\cdot)$ at the point $S\equiv 0$.
 \end{remark}

One check directly  that the function \eqref{sec:I.2}
satisfies the  conditions $\H_\zs{1})$-$\H_\zs{4})$. Another functions satisfying these conditions are given 
in Galtchouk, Pergamenshchikov , 2008.

\medskip 
\medskip
\vspace{5mm}
\section{ Main results}\label{sec:Ma}

Denote by $\cP_\zs{n}$  the family of distributions $p$ in $\bbr^n$ of the vectors
$(\xi_\zs{1},\ldots,\xi_\zs{n})'$
in the model \eqref{sec:I.1} such that the components $\xi_\zs{j}$ are jointly independent,
centered with unit variance and
\begin{equation}\label{sec:Ma.0}
\max_\zs{1\le k\le n}\,
\E\, \xi_\zs{k}^4\le l^*_\zs{n}\,,
\end{equation}
where $l^*_\zs{n}\ge 3$ is slowly increasing sequence, i.e.
it satisfies the condition \eqref{sec:Ad.6-1}.
It is easy to see that, for any $n\ge 1$, the centered gaussian distribution in $\bbr^n$ with
unit covariation matrix belongs to the family $\cP_\zs{n}$. We will denote by $q$ this gaussian distribution.

For any estimator $\wh{S}$, we define the following quadratic risk
\begin{equation}\label{sec:Ma.1}
\cR_\zs{n}(\wh{S},S)=\sup_\zs{p\in\cP_\zs{n}}\E_\zs{S,p}\|\wh{S}-S\|^2_\zs{n}\,,
\end{equation}
where $\E_\zs{S,p}$ is the expectation with respect to the distribution $\P_\zs{S,p}$
of the observations $(y_\zs{1},\ldots,y_\zs{n})$ with the fixed function $S$ and
the fixed  distribution $p\in\cP_\zs{n}$ of random variables $(\xi_\zs{j})_\zs{1\le j\le n}$ in the model \eqref{sec:I.1}.

Moreover, to make the risk independent of the design, we will make use of the risk with respect to the usual norm in $\cL_\zs{2}[0,1]$ \eqref{sec:Co.2} too,
i.e.
\begin{equation}\label{sec:Ma.1-1}
\cT_\zs{n}(\wh{S},S)=\sup_\zs{p\in\cP_\zs{n}}\E_\zs{S,p}\|\wh{S}-S\|^2\,.
\end{equation}

If an estimator $\wh{S}$ is defined only at the design points $(x_\zs{j})_\zs{1\le j\le n}$, then we extend it as step function onto the interval $[0,1]$ by setting $\wh{S}(x)=T(\wh{S}(x))$,
for all $0\le x\le 1$, where
\begin{equation}\label{sec:Ma.1-2}
T(f)(x)=f(x_\zs{1})\Chi_\zs{[0,x_{1}]}(x)
+\sum_{k=2}^n\,f(x_k)\Chi_\zs{(x_{k-1},x_k]}(x)\,.
\end{equation}

In Galtchouk, Pergamenshchikov, 2007, 2009
%\cite{GaPe1} 
the following non-asymptotic oracle inequality has been proved for the procedure 
\eqref{sec:Ad.7} .
\begin{theorem}\label{Th.sec:Ma.1}
Assume that $S\in W_\zs{r}^{1}$ for some $r>0$. Let $\wh{S}_\zs{*}$ be from \eqref{sec:Ad.7}.
Then, for any odd $n\ge 3$, the following oracle inequality holds:
\begin{equation}\label{sec:Ma.2}
\cR_\zs{n}(\wh{S}_\zs{*},S)\,
\le\,
\frac{1+3\rho-2\rho^2}{1-3\rho}
\,\min_\zs{\lambda\in\Lambda}\,
\cR_\zs{n}(\wh{S}_\zs{\lambda},S)\,
+\,\frac{1}{n}\,\cB_\zs{n}(\rho)\,,
\end{equation}
where $\rho=1/(3+L_n), L_n$ is from \eqref{sec:Ad.6-1}, the function $\cB_\zs{n}(\rho)$ is such that, for any $\nu>0$,
\begin{equation}\label{sec:Ma.3}
\lim_\zs{n\to\infty}\,
\frac{\cB_\zs{n}(\rho)}{n^\nu}\,=\,0\,.
\end{equation}
\end{theorem}

\begin{remark}\label{Re.sec:Ma.1}
Note that in Galtchouk, Pergamenshchikov, 2007, 2009, the oracle inequality is proved
for the model \eqref{sec:I.1}, where the random variables $(\xi_\zs{j})_\zs{1\le j\le n}$
are independent identically distributed. In fact, the result is true for 
independent random variables which are not identically distributed, i.e. for any
distribution of the random vector $(\xi_\zs{1},\ldots,\xi_\zs{n})'$  from $\cP_\zs{n}$.
\end{remark}

Now we formulate the main asymptotic results. To this end, 
for any function  $S\in W^k_\zs{r}$, we set
\begin{equation}\label{sec:Ma.4}
\gamma_k(S)\,=\,\Gamma^*_k\,r^{-(2k+1)}\,(\varsigma(S))^{2k/(2k+1)}\,,
\end{equation}
where
$$
\Gamma^*_k=(2k+1)^{-(2k+1)}\left(k/(\pi\,(k+1))\right)^{2k/(2k+1)}\,.
$$
 It is well known (see e.g., Nussbaum, 1985)
that the optimal rate of convergence is $n^{2k/(2k+1)}$ when 
the risk is taken uniformly over $W^k_\zs{r}$.

\begin{theorem}\label{Th.sec:Ma.2}
Assume that in the model \eqref{sec:I.1} the sequence
$(\sigma_\zs{j}(S))$ satisfies the condition 
$\H_1)$. Then, for any integer $k\ge 1$ and $r>0$, the estimator $\wh{S}_\zs{*}$ from \eqref{sec:Ad.7}
 satisfies the inequalities
\begin{equation}\label{sec:Ma.5}
\limsup_{n\to\infty}\,
n^{\frac{2k}{2k+1}}\,
\sup_\zs{S\in W^k_\zs{r}}\,
\frac{\cR_\zs{n}(\wh{S}_\zs{*},S)}{\gamma_k(S)}\,\le 1
\end{equation}
and
\begin{equation}\label{sec:Ma.5-1}
\limsup_{n\to\infty}\,
n^{\frac{2k}{2k+1}}\,
\sup_\zs{S\in W^k_\zs{r}}\,
\frac{\cT_\zs{n}(\wh{S}_\zs{*},S)}{\gamma_k(S)}\,\le 1\,.
\end{equation}
\end{theorem}
%The proof of this theorem is given in the next section.

The following result  gives the sharp lower bound for the risk \eqref{sec:Ma.1} and show that
$\gamma_k(S)$ is the \mbox{\it Pinsker} \mbox{\rm constant}.
\begin{theorem}\label{Th.sec:Ma.3}
Assume that in the model \eqref{sec:I.1}  the sequence
$(\sigma_\zs{j}(S))$
satisfies the conditions $\H_2)$- $\H_4)$. Then, for any integer $k\ge 1$ and $r>0$, the risks
\eqref{sec:Ma.1} and \eqref{sec:Ma.1-1}
admit the following asymptotic lower bounds, respectively,
\begin{equation}\label{sec:Ma.6}
\liminf_{n\to\infty}\,n^{\frac{2k}{2k+1}}\,\inf_{\wh{S}_\zs{n}}\,\sup_\zs{S\in W^k_\zs{r}}\,
\frac{\cR_\zs{n}(\wh{S}_\zs{n},S)}{\gamma_k(S)}\,
\ge 1
\end{equation}
and
\begin{equation}\label{sec:Ma.6-1}
\liminf_{n\to\infty}\,n^{\frac{2k}{2k+1}}\,\inf_{\wh{S}_\zs{n}}\,\sup_\zs{S\in W^k_\zs{r}}\,
\frac{\cT_\zs{n}(\wh{S}_\zs{n},S)}{\gamma_k(S)}\,
\ge 1\,.
\end{equation}
\end{theorem}

\begin{remark}\label{Re.sec:Ma.3}
To obtain the non-asymptotic oracle inequality \eqref{sec:Ma.2}, it is not necessary to make use of equidistant design points and the trigonometric basis. One may take any design points (deterministic or random) and any orthonormal basis satisfying \eqref{sec:Ad.0-1}. But
to obtain the property \eqref{sec:Ma.3} one needs to impose some technical conditions (see
Galtchouk, Pergamenshchikov, 2009).

Note that the results of Theorem~\ref{Th.sec:Ma.2} and Theorem~\ref{Th.sec:Ma.3} are based on equidistant design points and the trigonometric basis.
\end{remark}

\medskip
\section{ Upper bound}\label{sec:Up}

In this section we prove Theorem~\ref{Th.sec:Ma.2}.  To this end we will make use of the oracle inequality \eqref{sec:Ma.2}. We have to find an estimator from the family
\eqref{sec:Ad.2}-\eqref{sec:Ad.3} for which we can prove the upper bound \eqref{sec:Ma.5}. We start with the construction
of  such an estimator. First we put
\begin{equation}\label{sec:Up.0}
\wt{l}_\zs{n}=\inf\{i\ge 1\,:\,i\ve\ge \ov{r}(S)\}\wedge m
\quad\mbox{and}\quad
\ov{r}(S)=r/\varsigma(S)\,,
\end{equation}
where $a\wedge b=\min(a,b)$.
Then we choose an index from the set $\cA$ as
\begin{equation}\label{sec:Up.0-1}
\wt{\alpha}=(k,\wt{t}_\zs{n})\,,
\end{equation}
where $k$ is the parameter of the set $W^k_\zs{r}$ and $\wt{t}_\zs{n}=\wt{l}_\zs{n}\ve$.
Finally, we set
\begin{equation}\label{sec:Up.0-2}
\wt{S}=\wh{S}_\zs{\wt{\lambda}}
\quad\mbox{and}\quad
\wt{\lambda}=\lambda_\zs{\wt{\alpha}}\,.
\end{equation}
Now we formulate the upper bound \eqref{sec:Ma.5} for this estimator.
\begin{theorem}\label{Th.sec:Up.1} Assume that the condition $\H_1)$ holds. Then, for any integer $k\ge 1$ and $r>0$, 
\begin{equation}\label{sec:Up.1}
\limsup_{n\to\infty}\,n^{\frac{2k}{2k+1}}\,
\sup_{S\in W^k_r}\,
\frac{\cR_\zs{n}(\wt{S},S)}{\gamma_k(S)}\,\le 1\,.
\end{equation}
\end{theorem}
%The proof of this theorem is given in Section~\ref{sec:Pr}.
\begin{remark}\label{Re.sec:Up.1}
Note that the estimator
$\wt{S}$ belongs to the family \eqref{sec:Ad.2}-\eqref{sec:Ad.3}, but we can not use directly
 this estimator because the parameters $k$, $r$ and $\ov{r}(S)$ are unknown. We can use
this upper bound only via the oracle inequality \eqref{sec:Ma.2} proved for the
procedure \eqref{sec:Ad.7}.
\end{remark}

Now Theorem~\ref{Th.sec:Ma.1} and Theorem~\ref{Th.sec:Up.1} imply the upper bound
\eqref{sec:Ma.5}. To obtain the upper bound \eqref{sec:Ma.5-1} one needs the following auxiliary result.

\begin{lemma}\label{Le.sec:Up.1}
For any $0<\delta<1$ and  any estimator $\wh{S}_\zs{n}$ of $S\in W^k_\zs{r}$,
$$
\|\wh{S}_\zs{n}-S\|_\zs{n}^2\ge(1-\delta)\|T_\zs{n}(\wh{S})-S\|^2\,
-\,(\delta^{-1}\,-\,1)\,r/n^2\,,
$$
where the function
$T_\zs{n}(\wh{S})(\cdot)$ is defined in \eqref{sec:Ma.1-2}.
\end{lemma}
\noindent Proof of this lemma is given in
Galtchouk, Pergamenshchikov, 2008 ( Appendix A.1).

Now the inequality \eqref{sec:Ma.5} and this lemma imply the upper bound \eqref{sec:Ma.5-1}.
Hence Theorem~\ref{Th.sec:Ma.2}.

\medskip
\section{ Lower bound}\label{sec:Lo}

In this section we give the main steps of proving the lower bounds \eqref{sec:Ma.6} and
\eqref{sec:Ma.6-1}. We follow the commonly used scheme (see e.g. Nussbaum, 1985).
We begin with minorizing the minimax risk by a bayesian one constructed on a parametric functional family introduced in Section~\ref{subsec:Fa} ( see \eqref{Fa.4}) and using the
prior distribution \eqref{Fa.5}.
Further, a special modification of the van Trees inequality (see, Theorem~\ref{Th.sec:Tr.1}) yields
a lower bound for the bayesian risk depending on the chosen prior distribution, of course.
 Finally, in Section~\ref{subsec:Br}, we choose parameters of the prior distribution
(see \eqref{Fa.5}) providing the maximal value of the lower bound for the bayesian risk. This
value coincides with the Pinsker constant as it is shown in Section~\ref{subsec:Th.M.3}.
We emphasize that, by making use of the bayesian risk, one passes from the nonparametric regression model to a parametric one, for which the van Trees inequality holds. Note that this inequality is an extension of the Cramer-Rao inequality and gives a lower bound for the bayesian risk.

\subsection{ Lower bound for parametric heteroscedastic regression models} \label{subsec:Tr}

Let $(\bbr^n,\cB(\bbr^n),\P_\zs{\vartheta},\vartheta\in\Theta\subseteq\bbr^l)$ 
be the 
 statistical model relative to the observations $(y_\zs{j})_\zs{1\le j\le n}$ 
governed by the regression equation
\begin{equation}\label{subsec:Tr.1}
y_\zs{j}\,=\,S_\zs{\vartheta}(x_\zs{j})\,+\,\sigma_\zs{j}(\vartheta)\,\xi_\zs{j}\,,
\end{equation}
where $\xi_1,\ldots,\xi_\zs{n}$ are i.i.d. $\cN(0,1)$ random variables,
 $\vartheta=(\vartheta_1,\ldots,\vartheta_l)^\prime$ is a unknown parameter vector,
 $S_\zs{\vartheta}(x)$ is a unknown (or known) function 
and $\sigma_\zs{j}(\vartheta)=g(x_\zs{j},S_\vartheta)$, with the function
$g(x,S)$ defined in the condition $\H_\zs{1})$. Assume that a prior distribution
$\mu_\zs{\vartheta}$ of the parameter $\vartheta$ in $\bbr^l$ is defined by the density
 $\Phi(\cdot)$ of
the following form
$$
\Phi(z)\,=\,
\Phi(z_1,\ldots,z_l)=\prod_{i=1}^l\varphi_\zs{i}(z_\zs{i})\,,
$$
where $\varphi_\zs{i}$ is a continuously differentiable bounded density on $\bbr$ with
$$
I_\zs{i}=\int_{\bbr}\frac{\dot{\varphi}^2_\zs{i}(u)}{\varphi_\zs{i}(u)}\d u
<\infty\,.
$$
Let $\tau(\cdot)$ be a continuously differentiable $\bbr^l\to \bbr$ function
such that, for any $1\le i\le l$,
\begin{equation}\label{subsec:Tr.1-1}
\lim_\zs{|z_\zs{i}|\to\infty}\,\tau(z)\,\varphi_\zs{i}(z_\zs{i})=0
\quad\mbox{and}\quad
\int_{\bbr^l}\,\left|
\tau^{\prime}_\zs{i}(z)
\right|\,
\Phi(z)\d z<\infty\,,
\end{equation}
where 
$$ 
\tau^{\prime}_\zs{i}(z)=
(\partial/\partial z_\zs{i})\,\tau(z)\,.
$$
Let $\wh{\tau}_\zs{n}$ be an estimator of $\tau(\vartheta)$ based on
 observations $(y_\zs{j})_\zs{1\le j\le n}$.
 For any $\cB(\bbr^n\times\bbr^l)$ -
measurable integrable function 
$G(x,z), x\in\bbr^n, z\in\bbr^l$, 
we set
$$
\wt{\E}\,G(Y,\vartheta)\,=\,\int_{\bbr^l}\,\E_\zs{z}\,G(Y,z)\,
\Phi(z)\,\d z\,,
$$
where $\E_\zs{\vartheta}$ is the expectation with respect to the distribution 
$\P_\zs{\vartheta}$ of the vector $Y=(y_\zs{1},\ldots,y_\zs{n})$.
 Note that in this case
$$
\E_\zs{\vartheta}\,G(Y,\vartheta)\,=\,\int_\zs{\bbr^n}\,G(v,\vartheta)\,
f(v,\vartheta)\,\d v\,,
$$
where
\begin{equation}\label{subsec:Tr.2}
f(v,z)=\prod^n_\zs{j=1}\,\frac{1}{\sqrt{2\pi} \sigma_\zs{j}(z)}
\exp\left\{-\,\frac{(v_\zs{j}-S_\zs{z}(x_\zs{j}))^2}
{2\sigma_\zs{j}^2(z)}\right\}\,.
\end{equation}

We prove the following result.
\begin{theorem}\label{Th.sec:Tr.1}
Assume that the conditions $\H_1)-\H_2)$ hold.
Moreover, assume that 
 the function $S_\zs{z}(\cdot)$ with $z=(z_\zs{1},\ldots,z_\zs{l})'$ is
uniformly over $0\le x\le 1$
 differentiable with respect to
$z_\zs{i},\ 1\le i\le l$, i.e. 
for any $1\le i\le l$,
there exists
a function $S^{\prime}_\zs{z,i}\in \cC[0,1]$
such that
\begin{equation}\label{subsec:Tr.2-1}
\lim_\zs{h\to 0}\,
\max_\zs{0\le x\le 1}
\left|\left(
S_\zs{z+h\e_\zs{i}}(x)-S_\zs{z}(x)-S^{\prime}_\zs{z,i}(x)h\right)/h
\right|\,=0\,,
\end{equation}
where $\e_\zs{i}=(0,....,1,...,0)'$,  all coordinates are $0$, except the i-th equals to $1$ .
Then for any square integrable estimator $\wh{\tau}_\zs{n}$ of 
$\tau(\vartheta)$ and any $1\le i\le l$,
\begin{equation}\label{subsec:Tr.3}
\wt{\E}\left(\wh{\tau}_\zs{n}-\tau(\vartheta)\right)^2\ge
\frac{\ov{\tau}_\zs{i}^2}{F_\zs{i}\,+\,B_\zs{i}\,+\,I_\zs{i}}\,,
\end{equation}
where $\ov{\tau}_\zs{i}=
\int_{\bbr^l}\,
\tau^{\prime}_\zs{i}(z)\,
\Phi(z)\d z$, $F_\zs{i}=\sum_{j=1}^n
\int_\zs{\bbr^l}
\,(S^{\prime}_\zs{z,i}(x_\zs{j})/\sigma_\zs{j}(z))^2\,
\Phi(z)\d z$
and
$$
 B_\zs{i}=\,\frac{1}{2}\,\sum^n_\zs{j=1}\int_\zs{\bbr^l}
\frac{\wt{\L}^2_\zs{i}(x_\zs{j},S_\zs{z})}{\sigma^4_\zs{j}(S_\zs{z})}
\Phi(z)\d z\,, 
$$
$\wt{\L}_\zs{i}(x,z)=\L_\zs{x,S_\zs{z}}(S^{\prime}_\zs{z,i})$,
the operator $\L_\zs{x_,S}$ is defined in 
the condition $\H_2)$.
\end{theorem}

\begin{remark}\label{Re.sec:Tr.1}
Note that the inequality \eqref{subsec:Tr.3} is  some modification of the van Trees
inequality (see, Gill, Levit, 1995) adapted to the model \eqref{subsec:Tr.1}.
\end{remark}
\medskip

\subsection{ Parametric family of kernel functions}\label{subsec:Fa}

In this section we define and study some special parametric family of kernel function  
which will be used to prove the sharp lower bound \eqref{sec:Ma.6}.

Let us begin by kernel functions. We fix $\eta >0$ and we set 
\begin{equation}\label{Fa.1}
\chi_\zs{\eta}(x)
=\eta^{-1}\int_{\bbr}\Chi_\zs{(|u|\le 1-\eta)}\,V\left(\frac{u-x}{\eta}\right)\,\d u\,,
\end{equation}
where $\Chi_A$ is the indicator of a set $A$, the kernel $V\in\cC^{\infty}(\bbr)$  is 
such that
$$
V(u)=0
\quad\mbox{for}\quad |u|\ge 1
\quad\mbox{and}\quad 
\int^1_{-1}\,V(u)\,\d u=1\,.
$$
It is easy to see that the function $\chi_\zs{\eta}(x)$ possesses the properties :
\begin{align*}
&0\le \chi_\zs{\eta}\le 1\,,\quad \chi_\zs{\eta}(x)=1
\quad\mbox{for}\quad |x|\le 1-2\eta \quad\mbox{and}\\ 
&\chi_\zs{\eta}(x)=0
\quad\mbox{for}\quad |x|\ge 1\,.
\end{align*}
 Moreover,
for any $c>0$ and $\nu\ge 0$
\begin{equation}\label{Fa.2}
\lim_\zs{\eta\to 0}\,\sup_\zs{f\,:\,\|f\|_\zs{\infty}\le c}\,
\left|
\int_{\bbr}\,f(x)\,\chi_\zs{\eta}^\nu(x)\d x-\int_{-1}^{1}f(x) \d x
\right|
=0\,.
\end{equation}

We divide the interval $[0,1]$ into $M$ equal subintervals of length $2h$ and on each of them
we construct a  kernel-type function which equals to zero
at the boundary of the subinterval together with all derivatives. 
 This provides that the Fourier partial sums with respect
to the trigonometric basis in $\cL_\zs{2}[-1,1]$ give a natural parametric approximation to the
function on each subinterval. 
 
Let  $(e_\zs{j})_\zs{j\ge 1}$ be the trigonometric basis in $\cL_2[-1,1]$, i.e.
\begin{equation}\label{Fa.3}
e_\zs{1}=1/\sqrt{2}\,,\quad
e_\zs{j}(x)=\,Tr_\zs{j}\left(\pi [j/2] x\right)\,,\ j\ge 2\,,
\end{equation}
where the functions $(Tr_\zs{j})_\zs{j\ge 2}$ are defined in \eqref{sec:Ad.0-0}.

Now,
 for any array $z=\{(z_\zs{m,j})_\zs{1\le m\le M_\zs{n}\,,\,1\le j\le N_\zs{n}}\}$  we define
the following function
\begin{equation}\label{Fa.4}
S_\zs{z,n}(x)=\sum_{m=1}^{M_\zs{n}}\sum_{j=1}^{N_\zs{n}}\,z_\zs{m,j}\,D_\zs{m,j}(x)\,,
\end{equation}
where $D_\zs{m,j}(x)=e_\zs{j}\left(v_m(x)\right)\chi_\zs{\eta}\left(v_m(x)\right)$,
$$
v_m(x)=\frac{x-\wt{x}_m}{h_\zs{n}}\,,
\quad\wt{x}_m= 2mh_\zs{n}
\quad\mbox{and}\quad
M_\zs{n}=\left[1/(2h_\zs{n})\right]-1\,.
$$

We assume that the sequences 
$(N_\zs{n})_\zs{n\ge 1}$ and $(h_\zs{n})_\zs{n\ge 1}$,
satisfy the following conditions.

$\A_\zs{1})$
{\em The sequence $N_\zs{n}\to\infty$ as $n\to\infty$ and, for any $\nu>0$,
$$
\lim_\zs{n\to\infty}\,\frac{N^{\nu}_\zs{n}}{n}\,=\,0\,.
$$
Moreover, there exist $0<\delta_\zs{1}<1$ and $\delta_\zs{2}>0$
such that 
$$
h_\zs{n}=\,\O(n^{-\delta_\zs{1}})
\quad\mbox{and}\quad
h^{-1}_\zs{n}=\,\O(n^{\delta_\zs{2}})
\quad\mbox{as}\quad
n\to \infty\,.
$$
}
To define a prior distribution on the family of arrays,
we choose the following random array 
$\vartheta=\{(\vartheta_\zs{m,j})_\zs{ 1\le m\le M_\zs{n}\,,\, 1\le j\le N_\zs{n}}\}$ 
with
\begin{equation}\label{Fa.5}
\vartheta_\zs{m,j}\,=\,t_\zs{m,j}\,\zeta_\zs{m,j}\,,
\end{equation}
where  $(\zeta_\zs{m,j})$ are i.i.d. $\cN(0,1)$ random variables and 
$(t_\zs{m,j})$
are some nonrandom positive coefficients. We make use of gaussian variables since they
possess the minimal Fisher information and therefore maximize the lower bound
\eqref{subsec:Tr.3}.
 We set
\begin{equation}\label{Fa.6}
t^*_\zs{n}=\,\max_\zs{1\le m\le M_\zs{n}}
\sum^{N_\zs{n}}_\zs{j=1}\,t_\zs{m,j}\,.
\end{equation}
We assume that the coefficients $(t_\zs{m,j})_\zs{ 1\le m\le M_\zs{n}\,,\, 1\le j\le N_\zs{n}}$
satisfy the following conditions.

$\A_\zs{2})$
{\em There exists  a sequence of positive numbers $(d_\zs{n})_\zs{n\ge 1}$ such that
\begin{equation}\label{Fa.7}
\lim_\zs{n\to\infty}
\frac{d_\zs{n}}{h_\zs{n}^{2k-1}}\,
\sum^{M_\zs{n}}_\zs{m=1}\sum^{N_\zs{n}}_\zs{j=1}\,t^2_\zs{m,j}\,
j^{2(k-1)}=0\,,
\quad
\lim_\zs{n\to\infty}\,\sqrt{d_\zs{n}}\,t^*_\zs{n}=0\,,
\end{equation}
moreover, for any $\nu> 0$,
$$
\lim_\zs{n\to\infty}\,n^{\nu}\,\exp\{{-d_\zs{n}/2}\}=0\,.
$$
}

$\A_\zs{3})$
{\em For some $0<\epsilon<1$
$$
\limsup_\zs{n\to\infty}
\frac{1}{h_\zs{n}^{2k-1}}\,\sum^{M_\zs{n}}_\zs{m=1}\sum^{N_\zs{n}}_\zs{j=1}\,t^2_\zs{m,j}\,
j^{2k}\,
\le (1-\epsilon)r
\left(\frac{2}{\pi}\right)^{2k}
\,.
$$
}

$\A_\zs{4})$
{\em There exists $\epsilon_\zs{0}>0$ such that
$$
\lim_\zs{n\to\infty}
\frac{1}{h_\zs{n}^{4k-2+\epsilon_\zs{0}}}\,
\sum^{M_\zs{n}}_\zs{m=1}\sum^{N_\zs{n}}_\zs{j=1}\,t^4_\zs{m,j}\,j^{4k}\,=0\,.
$$
}

\begin{prop}\label{sec:Pr.Fa.1}
Let the conditions $\A_\zs{1})$--$\A_\zs{2})$. Then,
for any $\nu>0$ and for any $\delta>0$,
$$
\lim_\zs{n\to\infty}\,n^\nu\,\max_\zs{0\le l\le k-1}
\P\left(\|S^{(l)}_\zs{\vartheta,n}\|>\delta\right)=0\,.
$$
\end{prop}

\begin{prop}\label{sec:Pr.Fa.2}
Let the conditions
$\A_\zs{1})$--$\A_\zs{4})$. Then, for any $\nu>0$,
$$
\lim_\zs{n\to\infty}\,n^\nu\,
\P(S_\zs{\vartheta,n}\notin W^{k}_\zs{r})\,=0\,.
$$
\end{prop}

\medskip

\begin{prop}\label{sec:Pr.Fa.3}
Let the conditions
$\A_\zs{1})$--$\A_\zs{4})$. Then, for any $\nu>0$,
$$
\lim_\zs{n\to\infty}\,n^\nu\,
\E\,\|S_\zs{\vartheta,n}\|^2\,
\left(
\Chi_\zs{\{S_\zs{\vartheta,n}\notin W^{k}_\zs{r}\}}\,+\,
\Chi_\zs{\Xi^c_\zs{n}}
\right)
=0\,.
$$
\end{prop}
\medskip

\begin{prop}\label{sec:Pr.Fa.4}
Let the conditions
$\A_\zs{1})$--$\A_\zs{4})$. Then, for any function $g$
satisfying the conditions \eqref{sec:Co.6} and $\H_\zs{4})$, 
$$
\lim_\zs{n\to\infty}\,\sup_\zs{0\le x\le 1}\,
\E\,\left|\,g^{-2}(x,S_\zs{\vartheta,n})-g^{-2}_\zs{0}(x)\right|=0\,.
$$
\end{prop}

Proofs of Propositions~\ref{sec:Pr.Fa.1}--\ref{sec:Pr.Fa.4}
are given in Galtchouk, Pergamenshchikov, 2008 (Appendix A.2--A.6).
\medskip

\subsection{Bayes risk}\label{subsec:Br}

Now we will obtain the lower bound for the bayesian risk that yields the lower bound
\eqref{sec:Ma.6-1} for the minimax risk.
% V etom razdele my poluchim nijnuyu ocenku dlya baesovskogo riska, kotoryi ispol'zuetsy
%v poluchenii nujnei granitsy \eqref{sec:Ma.6-1}.

We make use of the sequence of random functions
$(S_\zs{\vartheta,n})_\zs{n\ge 1}$ defined in \eqref{Fa.4}-\eqref{Fa.5}
with the coefficients $(t_\zs{m,j})$ satisfying the conditions $\A_\zs{1})$--$\A_\zs{4})$
which will be chosen later.

For any estimator $\wh{S}_\zs{n}$ we introduce now the corresponding Bayes risk
\begin{equation}\label{sec:Lo.3}
\cE_\zs{n}(\wh{S}_\zs{n})=
\int_\zs{\bbr^l}\,\E_\zs{S_\zs{z,n},q}\|\wh{S}_\zs{n}-S_\zs{z,n}\|^2\,\mu_\zs{\vartheta}(\d z)\,,
\end{equation}
where the kernel family $(S_\zs{z,n})$ is defined in \eqref{Fa.4}, $\mu_\zs{\vartheta}$ denotes 
 the distribution of the random array
 $\vartheta$ defined by \eqref{Fa.5} in $\bbr^l$ with $l= M_\zs{n}N_\zs{n}$.

We remember that $q$ is a centered gaussian distribution in $\bbr^n$ with unit covariation
 matrix.

First of all, we replace the functions
$\wh{S}_\zs{n}$ and $S$ by their Fourier series with respect to the basis
$$
\wt{e}_\zs{m,i}(x)=(1/\sqrt{h})\,e_\zs{i}\left(v_m(x)\right)\,
\Chi_\zs{\left(|v_m(x)|\le 1\right)}\,.
$$
By making use of this basis
we can estimate the norm $\|\wh{S}_\zs{n}-S_\zs{z,n}\|^2$ from below as
$$
\|\wh{S}_\zs{n}-S_\zs{z,n}\|^2 \ge
 \sum_{m=1}^{M_\zs{n}}\sum_{j=1}^{N_\zs{n}}\,(\wh{\tau}_\zs{m,j}\,-\,
\tau_\zs{m,j}(z))^2\,,
$$
where
$$
\wh{\tau}_\zs{m,j}=\int_0^1\,\wh{S}_\zs{n}(x)\wt{e}_\zs{m,j}(x)\d x
\quad\mbox{and}\quad 
\tau_\zs{m,j}(z)=\int_0^1\,S_\zs{z,n}(x)\wt{e}_\zs{m,j}(x)\,\d x\,.
$$
Moreover, from the definition \eqref{Fa.4}  one gets
$$
\tau_\zs{m,j}(z)
=\sqrt{h}\sum_{i=1}^{N_\zs{n}}\,z_\zs{m,i}\int_{-1}^1\,e_\zs{i}(u)e_\zs{j}(u)\chi_\zs{\eta}(u)\,\d u\,.
$$
It is easy to see that the functions $\tau_\zs{m,j}(\cdot)$
satisfy the condition \eqref{subsec:Tr.1-1} for gaussian prior densities. In this case
(see the definition in \eqref{subsec:Tr.3}) we have
$$
\ov{\tau}_\zs{m,j}=
(\partial/\partial z_\zs{m,j}) \tau_\zs{m,j}(z)
=\sqrt{h} \ov{e}_\zs{j}(\chi_\zs{\eta})\,,
$$
where 
\begin{equation}\label{sec:Lo.6}
\ov{e}_\zs{j}(f)=\int_{-1}^1\,e^2_\zs{j}(v)\,f(v)\,\d v\,.
\end{equation}
Now to obtain a lower bound for
the Bayes risk 
$\cE_\zs{n}(\wh{S}_\zs{n})$
 we make use of Theorem~\ref{Th.sec:Tr.1} which implies that
\begin{equation}\label{sec:Lo.7}
\cE_\zs{n}(\wh{S}_\zs{n})\,\ge
\sum_{m=1}^{M_\zs{n}}\sum_{j=1}^{N_\zs{n}}\,
\frac{h \ov{e}^2_\zs{j}(\chi_\zs{\eta})}{F_\zs{m,j}+B_\zs{m,j}\,  
+t^{-2}_\zs{m,j}}\,,
\end{equation}
where $F_\zs{m,j}=\sum_{i=1}^n\,D^2_\zs{m,j}(x_\zs{i})\,
\E\,g^{-2}(x_\zs{i},S_\zs{\vartheta,n})$ and
$$
B_\zs{m,j}=
\frac{1}{2}\sum_{i=1}^n\,
\E\,
\frac{\wt{\L}^2_\zs{m,j}(x_\zs{i},S_\zs{\vartheta,n})}{g^4(x_\zs{i},S_\zs{\vartheta,n})}
\quad\mbox{with}\quad
\wt{\L}_\zs{m,j}(x,S)=\L_\zs{x,S} \left(D_\zs{m,j}\right)\,.
$$

By making use of Proposition~\ref{sec:Pr.Fa.4} we can show that
\begin{equation}\label{sec:Lo.8}
\lim_\zs{n\to\infty}\,
\sup_\zs{1\le m\le M_\zs{n}}\,\sup_\zs{1\le j\le N_\zs{n}}
\left|
\frac{1}{nh}\,
F_\zs{m,j}\,-\,\ov{e}_\zs{j}(\chi_\zs{\eta}^2)\,
g^{-2}_\zs{0}(\wt{x}_m)
\right|=0
\end{equation}
and
\begin{equation}\label{sec:Lo.9}
\lim_\zs{n\to\infty}
\frac{1}{nh}
\sup_\zs{1\le m\le M_\zs{n}}\,\sup_\zs{1\le j\le N_\zs{n}}\,B_\zs{m,j}\,=\,0\,.
\end{equation}
The detailed proof of these equalities is given in Galtchouk, Pergamenshchikov, 2008 (Appendix A.8).

This means that, for any $\nu>0$
and for sufficiently large $n$,
$$
\sup_\zs{1\le m\le M_\zs{n}}\,\sup_\zs{1\le j\le N_\zs{n}}
\frac{F_\zs{m,j}+B_\zs{m,j}+t^{-2}_\zs{m,j}}
{nh \ov{e}_\zs{j}(\chi_\zs{\eta}^2) g^{-2}_\zs{0}(\wt{x}_m)
+t^{-2}_\zs{m,j}} \le 1+\nu\,,
$$
where $\wt{x}_\zs{m}$ is defined in \eqref{Fa.4}. Therefore, if we  denote in \eqref{sec:Lo.7} 
$$
\kappa^2_\zs{m,j}=\,nh\,g^{-2}_\zs{0}(\wt{x}_m)\,t^2_\zs{m,j}
\quad\mbox{and}\quad
\psi_\zs{j}(\eta,y)
=\,
\frac{\ov{e}^2_\zs{j}(\chi_\zs{\eta})y}{\ov{e}_\zs{j}^2(\chi_\zs{\eta}^2)y+1}
$$
we obtain, for sufficiently large $n$,
$$
n^{\frac{2k}{2k+1}}\cE_\zs{n}(\wh{S}_\zs{n})\,\ge
\frac{n^{-\frac{1}{2k+1}}}{1+\nu}
\,
\sum_{m=1}^{M_\zs{n}}\,g^2_\zs{0}(\wt{x}_m)
\,
\sum_{j=1}^{N_\zs{n}}\,\psi_\zs{j}(\eta,\kappa^2_\zs{m,j})\,.
$$
Moreover, the  property
\eqref{Fa.2} implies
\begin{equation}\label{sec:Lo.10}
\lim_\zs{\eta\to 0}\,\sup_\zs{N\ge 1}\sup_\zs{(y_\zs{1},\ldots,y_\zs{N})\in \bbr^{N}_\zs{+}}
\left|
\frac{\sum_{j=1}^{N}\,\psi_\zs{j}(\eta,y_\zs{j})}{\Psi_\zs{N}(y_\zs{1},\ldots,y_\zs{N})}
\,-\,1
\right|\,=\,0\,,
\end{equation}
where
$$
\Psi_\zs{N}(y_\zs{1},\ldots,y_\zs{N})=\sum^N_\zs{j=1}\,\frac{y_\zs{j}}{y_\zs{j}+1}\,.
$$
Therefore we can write that, for sufficiently large $n$,
\begin{equation}\label{sec:Lo.11}
n^{\frac{2k}{2k+1}}\cE_\zs{n}(\wh{S}_\zs{n})\,\ge
\frac{1-\nu}{1+\nu}\,n^{-\frac{1}{2k+1}}
\,
\sum_{m=1}^{M_\zs{n}}\,g^2_\zs{0}(\wt{x}_m)
\,\Psi_\zs{N_\zs{n}}(\kappa^2_\zs{m,1},\ldots,\kappa^2_\zs{m,N_\zs{n}})\,.
\end{equation}

Obviously, to obtain a "good" lower bound for the risk 
$\cE_\zs{n}(\wh{S}_\zs{n})$
 one needs to maximize the right-hand side of the inequality \eqref{sec:Lo.11}. Hence we
choose the coefficients $(\kappa^2_\zs{m,j})$ by maximizing  the function
$\Psi_\zs{N}$, i.e.
$$
\max_\zs{y_\zs{1},\ldots,y_\zs{N}}\,\Psi_\zs{N}(y_\zs{1},\ldots,y_\zs{N})
\quad\mbox{subject to}\quad
\sum^N_\zs{j=1}y_\zs{j}j^{2k}\le R\,.
$$
The parameter $R>0$ will be chosen later to satisfy the condition $\A_\zs{3})$. 
By the Lagrange multipliers method it is easy to find that
the solution of this problem is given by
\begin{equation}\label{sec:Lo.12}
y^*_\zs{j}(R)=a^*(R)\,j^{-k}-1
\end{equation}
with
$$
a^*(R)=\frac{1}{\sum^N_\zs{i=1}\,i^{k}}
\left(R+\sum^N_\zs{i=1}i^{2k}\right)
\quad\mbox{and}\quad
1\le j\le N\,.
$$

To obtain a positive solution in \eqref{sec:Lo.12} we need to impose  the following condition
\begin{equation}\label{sec:Lo.13}
R> \,N^{k}\,\sum^N_\zs{i=1} i^{k}-\sum^N_\zs{i=1}i^{2k}\,.
\end{equation}
Moreover, from the condition $\A_\zs{3})$ we  obtain that
\begin{equation}\label{sec:Lo.14}
R\le \frac{2^{2k+1}(1-\ve)r\,
n\,h^{2k+1}_\zs{n}}{\pi^{2k}\wh{g}_\zs{0}}:=R^*_\zs{n}\,,
\end{equation}
where
$$
\wh{g}_\zs{0}=2h_\zs{n}\,\sum_{m=1}^{M_\zs{n}}\,g^2_\zs{0}(\wt{x}_m)\,.
$$
Note that, by the condition $\H_\zs{4})$, the function 
$g_\zs{0}(\cdot)=g(\cdot,S_\zs{0})$ is continuous on 
the interval $[0,1]$, therefore, 
\begin{equation}\label{sec:Lo.15}
\lim_\zs{n\to\infty}\wh{g}_\zs{0}=\int^1_\zs{0}g^2(x,S_\zs{0})\d x=\varsigma(S_\zs{0})
\quad\mbox{with}\quad
S_\zs{0}\equiv 0\,.
\end{equation}

Now we have to choose the sequence $(h_\zs{n})$. Note that if we put in \eqref{Fa.5}
\begin{equation}\label{sec:Lo.15-1}
t_\zs{m,j}=g_\zs{0}(\wt{x}_m)\sqrt{y^*_\zs{j}(R)}/\sqrt{nh_\zs{n}}
\quad\mbox{i.e.}\quad
\kappa�_\zs{m,j}\,=\,y^*_\zs{j}(R)\,,
\end{equation}
we can rewrite the inequality \eqref{sec:Lo.11} as
\begin{equation}\label{sec:Lo.16}
n^{\frac{2k}{2k+1}}\cE_\zs{n}(\wh{S}_\zs{n})\,\ge
\frac{(1-\nu)}{(1+\nu)}\,\frac{\wh{g}_\zs{0}\,\Psi^*_\zs{N_\zs{n}}(R)}{2h_\zs{n}}\, n^{-\frac{1}{2k+1}}\,,
\end{equation}
where
$$
\Psi^*_\zs{N}(R)=
N-\frac{\left(\sum^{N}_\zs{j=1}j^{k}\right)^2}{R+\sum^{N}_\zs{j=1}j^{2k}}\,.
$$
It is clear that
$$
k^2/(k+1)^2\le \liminf_\zs{N\to\infty}\inf_\zs{R>0}\Psi^*_\zs{N}(R)/N
\le \limsup_\zs{N\to\infty}\sup_\zs{R>0}\Psi^*_\zs{N}(R)/N\le 1\,.
$$
Therefore, to obtain a positive finite asymptotic lower bound in \eqref{sec:Lo.16}
we have to take the parameter $h_\zs{n}$ as
\begin{equation}\label{sec:Lo.17}
h_\zs{n}=h_\zs{*}n^{-1/(2k+1)}N_\zs{n}
\end{equation}
with some positive  coefficient $h_\zs{*}$. Moreover,  the conditions 
\eqref{sec:Lo.13}-\eqref{sec:Lo.14} imply that,
for sufficiently large $n$,
$$
(1-\ve)r\,\frac{2^{2k+1}}{\pi^{2k}}
\,\frac{1}{\wh{g}_\zs{0}}\,h^{2k+1}_\zs{*}
>
\frac{1}{N^{k+1}_\zs{n}}\sum^{N_\zs{n}}_\zs{j=1}j^k-\frac{1}{N^{2k+1}_\zs{n}}
\sum^{N_\zs{n}}_\zs{j=1}j^{2k}\,.
$$
Now taking into account that, for sufficiently large $n$,
$$
\wh{g}_\zs{0}\,
\frac{1}{N_\zs{n}}\sum^{N_\zs{n}}_\zs{j=1}
\left(
(j/N_\zs{n})^k
-
(j/N_\zs{n})^{2k}
\right)
<\,
\frac{(1+\ve)\,\varsigma(S_\zs{0})k}{(k+1)(2k+1)}\,,
$$
we obtain the following condition on $h_\zs{*}$:
\begin{equation}\label{sec:Lo.18}
h_\zs{*}\ge (\upsilon^*_\zs{\ve})^{1/(2k+1)}\,,
\end{equation}
where
$$
\upsilon^*_\zs{\ve}=\,\frac{(1+\ve)k}{c^*_\zs{\ve}(k+1)(2k+1)}
\quad\mbox{and}\quad
c^*_\zs{\ve}=\frac{2^{2k+1}(1-\ve)r}{\pi^{2k}\varsigma(S_\zs{0})}\,.
$$
To maximize the function $\Psi^*_\zs{N_\zs{n}}(R)$
on the right-hand side of the inequality \eqref{sec:Lo.16}, we take  $R=R^*_\zs{n}$ defined in \eqref{sec:Lo.14}.
Therefore, we obtain that
\begin{equation}\label{sec:Lo.19}
\liminf_\zs{n\to\infty}\,\inf_\zs{\wh{S}_\zs{n}}\,
n^{2k/(2k+1)}\cE_\zs{n}(\wh{S}_\zs{n})\,\ge
\,\varsigma(S_\zs{0})\,F(h_\zs{*})/2\,,
\end{equation}
where
$$
F(x)=\frac{1}{x}-\frac{2k+1}{(k+1)^2(c^*_\zs{\ve}(2k+1)x^{2k+2}+x)}\,.
$$
Furthermore, taking into account that
$$
F^{\prime}(x)=-
\frac{(c^*_\zs{\ve}(2k+1)(k+1)x^{2k+1}-k)^2}{(k+1)^2(c^*_\zs{\ve}(2k+1)x^{2k+2}+x)^2}\le 0\,,
$$
we get
$$
\max_\zs{h_\zs{*}\ge (\upsilon^*_\zs{\ve})^{1/(2k+1)}}F(h_\zs{*})=
F((\upsilon^*_\zs{\ve})^{1/(2k+1)})=
\frac{(1+\ve^\prime)k}{k+1}
\,(\upsilon^*_\zs{\ve})^{-1/(2k+1)}\,,
$$
where
\begin{equation}\label{sec:Lo.19-1}
\ve^\prime=\frac{\ve}{2k+\ve k+1}\,.
\end{equation}

This means that to obtain in \eqref{sec:Lo.19} the maximal lower bound, one has to take
in \eqref{sec:Lo.17}
\begin{equation}\label{sec:Lo.20}
h_\zs{*}=(\upsilon^*_\zs{\ve})^{1/(2k+1)}\,.
\end{equation}
It is important to note that if one defines the prior distribution $\mu_\zs{\vartheta}$
in the bayesian risk \eqref{sec:Lo.3} by formulas \eqref{Fa.5}, \eqref{sec:Lo.15-1}, \eqref{sec:Lo.17} and \eqref{sec:Lo.20}, then the bayesian risk should depend on a parameter
$0<\ve<1$, i.e. $\cE_\zs{n}=\cE_\zs{\ve,n}$.
Therefore, the inequality \eqref{sec:Lo.19} implies that, for any $0<\ve<1$,
\begin{equation}\label{sec:Lo.21}
\liminf_\zs{n\to\infty}\,\inf_\zs{\wh{S}_\zs{n}}\,
n^{2k/(2k+1)}\cE_\zs{\ve,n}(\wh{S}_\zs{n})
\ge
\frac{(1+\ve^\prime)(1-\ve)^{1/(2k+1)}}{(1+\ve)^{1/(2k+1)}}\,\gamma_k(S_\zs{0})\,,
\end{equation}
where the function $\gamma_k(S_\zs{0})$ is defined in \eqref{sec:Ma.4} for $S_\zs{0}\equiv 0$.

Now to end the definition of the sequence of the random functions
$(S_\zs{\vartheta,n})$ defined by \eqref{Fa.4} and 
\eqref{Fa.5}, one has to define the sequence 
$(N_\zs{n})$.
 Let us remember that we make use of the  sequence 
$(S_\zs{\vartheta,n})$ 
with the coefficients
$(t_\zs{m,j})$
constructed in \eqref{sec:Lo.15-1} for $R=R^*_\zs{n}$ given in 
\eqref{sec:Lo.14} and for the sequence $h_\zs{n}$ given by 
\eqref{sec:Lo.17} and \eqref{sec:Lo.20}, for some fixed $0<\ve<1$.

  We will  choose the sequence $(N_\zs{n})$ 
to satisfy the conditions $\A_\zs{1})$--$\A_\zs{4})$. One can take, for example, 
$N_\zs{n}=[\ln^{4} n]+1$. Then the condition $\A_\zs{1})$ is trivial. Moreover,
taking into account  that in
this case
$$
R^*_\zs{n}=\frac{2^{2k+1}(1-\ve)r}{\pi^{2k} \wh{g}_\zs{0}}\upsilon^*_\zs{\ve} N^{2k+1}_\zs{n}
=
\frac{\varsigma(S_\zs{0})}{\wh{g}_\zs{0}}
\frac{k}{(k+1)(2k+1)}\,N^{2k+1}_\zs{n}\,,
$$
we find, thanks to the convergence
\eqref{sec:Lo.15},
$$
\lim_\zs{n\to\infty}\,
\dfrac{R^*_\zs{n}+\sum^{N_\zs{n}}_\zs{j=1}j^{2k}}{N^{k}_\zs{n} \sum^{N_\zs{n}}_\zs{j=1}j^{k}}=1\,.
$$
Therefore, the solution \eqref{sec:Lo.12}, for sufficiently large $n$, 
satisfies the following inequality:
$$
\max_\zs{1\le j\le N_\zs{n}}
y^*_\zs{j}(R^*_\zs{n})\,j^k
\le 2 N^k_\zs{n}\,.
$$
Now it is easy to see that  the condition 
 $\A_\zs{2})$ holds with $d_\zs{n}=\sqrt{N_\zs{n}}$ and 
the condition
$\A_\zs{4})$ holds for arbitrary 
$0<\epsilon_\zs{0}<1$.
As to the condition $\A_\zs{3})$, note that, in view of the definition of $t_\zs{m,j}$ 
in \eqref{sec:Lo.15-1}, we get
\begin{align*}
\frac{1}{h_\zs{n}^{2k-1}}\,\sum^{M_\zs{n}}_\zs{m=1}\sum^{N_\zs{n}}_\zs{j=1}\,t^2_\zs{m,j}\,j^{2k}
&=\frac{1}{2 nh_\zs{n}^{2k+1}}\,\wh{g}_\zs{0}\,
\sum^{N_\zs{n}}_\zs{j=1}\,y^*_\zs{j}(R^*_\zs{n})\,j^{2k}\\
&
=\frac{R^*_\zs{n} \wh{g}_\zs{0}}{N^{2k+1}_\zs{n} 2\upsilon^*_\zs{\ve}}
=
(1-\epsilon)r
\left(\frac{2}{\pi}\right)^{2k}\,.
\end{align*}
Hence the condition $\A_\zs{3})$.

\medskip
\section{Estimation of non periodic function}\label{sec:Np}

Now we consider the estimation problem of a non periodic regression  function $S$ in the model
\eqref{sec:I.1}. In this case we will estimate the function $S$ on any interior interval $[a,b]$
of $[0,1]$, i.e. for $0<a<b<1$.
 
It should be pointed out that at the boundary points $x=0$ and $x=1$, one must to make use of
kernel estimators (see Brua, 2007).
 
Let now $\chi$ be a infinitely differentiable $[0,1]\to\bbr_\zs{+}$ function such that 
$\chi(x)=1$ for $a\le x\le b$ and $\chi^{(k)}(0)=\chi^{(k)}(1)=0$ for all $k\ge 0$, for example,
$$
\chi(x)=\frac{1}{\eta}\,\int^{\infty}_\zs{-\infty}\,
V\left(\frac{x-z}{\eta}\right)\,
\Chi_\zs{[a',b']}(z)\,
\d z\,,
$$
where $V$ is some kernel function introduced in \eqref{Fa.1},
$$
a'=\frac{a}{2}\,,\quad b'=\frac{b}{2}+\frac{1}{2}
\quad\mbox{and}\quad
\eta=\frac{1}{4}\min(a\,,\,1-b)\,.
$$
Multiplying the equation \eqref{sec:I.1} by the function $\chi(\cdot)$ and simulating some
i.i.d. $\cN(0,1)$ sequence $(\zeta_\zs{j})_\zs{1\le j\le n}$ one comes to the estimation
problem of the periodic regression function $\wt{S}(x)=S(x)\chi(x)$ in the model 
$$
\wt{y}_\zs{j}=\wt{S}(x_\zs{j})+\wt{\sigma}_\zs{j}(S)\,\wt{\xi}_\zs{j}\,, 
$$
where 
$\wt{\sigma}_\zs{j}(S)=\sqrt{\sigma^2_\zs{j}(S)+\epsilon^2}$,
$$
\wt{\xi}_\zs{j}=\frac{\sigma_\zs{j}(S)}{\wt{\sigma}_\zs{j}(S)}\,\xi_\zs{j}+
\frac{\epsilon}{\wt{\sigma}_\zs{j}(S)}\,\zeta_\zs{j}\,.
$$
and $\epsilon>0$ is some sufficiently small parameter.

 It is easy to see that if the sequence $(\sigma_\zs{j}(S))_\zs{1\le j\le n}$ satisfies the conditions 
$\H_\zs{1})-\H_\zs{4})$, then the sequence $(\wt{\sigma}_\zs{j}(S))_\zs{1\le j\le n}$
satisfies these conditions as well with 
$$
\wt{\sigma}_\zs{j}(S)=\wt{g}(x_\zs{j},S)=\sqrt{g^2(x_\zs{j},S)\chi^2(x_\zs{j})+\epsilon^2}\,.
$$
\medskip
\section{Conclusion}\label{sec:Cn}
In conclusion, it should be noted that this paper completes the investigation of the
estimation problem for the nonparametric regression function in the heteroscedastic regression
model \eqref{sec:I.1} in the case of quadratic risk. It is shown that the adaptive procedure
\eqref{sec:Ad.7} satisfies the non asymptotic oracle inequality and it is asymptotically
efficient for estimating a periodic regression function as well. 
 From practical point of view, the procedure \eqref{sec:Ad.7} gives an acceptable accuracy
even for small samples as it is shown via simulations by  Galtchouk, Pergamenshchikov, 2009.

\medskip
\section{ Proofs}\label{sec:Pr}
\subsection{ Properties of the trigonometric basis}\label{subsec:A.7}

\begin{lemma}\label{Le.A.1} For any function $S\in  W_r^k$,
\begin{equation}\label{A.0-3}
\sup_{n\ge 1}\sup_{1\le m\le n-1}\,m^{2k}\,
\left(\sum_{j=m+1}^{n}\,\theta_\zs{j,n}^2
\right)
\,\le\,\frac{4r}{\pi^{2(k-1)}}\,.
\end{equation}
\end{lemma}

\begin{lemma}\label{Le.A.2}
For any $m\ge 0$,
\begin{equation}\label{A.0-2}
\sup_{N\ge 2}\quad\sup_{x\in [0,1]}N^{-m}\left|\sum_{l=2}^{N}\,l^m
\,
\left( \phi^2_\zs{l}(x)-1 \right)\,
\right|\le 2^m\,.
\end{equation}
%where $\ov{\phi}_\zs{l}(x)=\phi^2_\zs{l}(x)-1$.
\end{lemma}
Proofs of Lemma~\ref{Le.A.1} and Lemma~\ref{Le.A.2} are given in Galtchouk, Pergamenshchikov, 2007.

\begin{lemma}\label{Le.A.3}
Let $\theta_\zs{j,n}$ and $\theta_\zs{j}$ be the Fourier coefficients defined in \eqref{sec:Ad.1}
and \eqref{sec:Co.3-1}, respectively.
 Then, for $1\le j\le n$
and $n\ge 2$,
\begin{equation}\label{A.0-1}
\sup_{S\in W^1_r}\,|\theta_\zs{j,n}-\theta_\zs{j}|\,\le\, 2\pi\,\sqrt{r}\,j/n\,.
\end{equation}
\end{lemma}
\noindent{\bf Proof.}
Indeed, we have
\begin{align*}
|\theta_\zs{j,n}-\theta_\zs{j}|&\,=\,
\left|\sum_{l=1}^n\int_{x_\zs{l-1}}^{x_\zs{l}}\,\left(S(x_l)\phi_\zs{j}(x_l)-S(x)\phi_\zs{j}(x)\right)\d x\right|\\
&\le\,n^{-1}\,\sum_{l=1}^n\int_{x_\zs{l-1}}^{x_\zs{l}}\,
\left(|\dot{S}(z)\phi_\zs{j}(z)|\,+ \,|S(z)\dot{\phi}_\zs{j}(z)|\right)\d z
\\
&=\,n^{-1}\,\int_{0}^{1}\,
\left(|\dot{S}(z)|\,|\phi_\zs{j}(z)|\,+ \,|S(z)|\,|\dot{\phi}_\zs{j}(z)|\right)\d z\,.
\end{align*}
By making use of the Bounyakovskii-Cauchy-Schwarz inequality we get
\begin{align*}
|\theta_\zs{j,n}-\theta_\zs{j}|&\le\,
n^{-1}\,\left(\|\dot{S}\|\,\|\phi\|\,+\,\|\dot{\phi}\|\,\|S\|\right)\\
&\le\,
n^{-1}\,\left(\|\dot{S}\|\,+\,\pi\,j\,\,\|S\|\right)\,.
\end{align*}
The definition of the class $W^1_r$  implies (\ref{A.0-1}).
Hence Lemma~\ref{Le.A.1}.
\endproof

\subsection{ Proof of Theorem~\ref{Th.sec:Up.1}}\label{subsec:Th.sec:Up.1}
To prove the theorem 
we will adapt to the  heteroscedastic case the corresponding proof
 from Nussbaum, 1985.

First, from \eqref{sec:Ad.2}  we obtain that, for any $p\in\cP_\zs{n}$,
\begin{equation}\label{sec:Up.1-1}
\E_\zs{S,p}\,\|\wt{S}-S\|_\zs{n}^2=
\sum_{j=1}^{n}\,(1\,-\,\wt{\lambda}_\zs{j})^2
\theta^2_\zs{j,n}
+
\frac{1}{n}
\sum_{j=1}^n\,\wt{\lambda}_\zs{j}^2\varsigma_\zs{j,n}(S)\,,
\end{equation}
where
 $$
\varsigma_\zs{j,n}(S)=\frac{1}{n}\,
\sum^n_\zs{l=1}\sigma^2_\zs{l}(S)\phi^2_\zs{j}(x_\zs{l})\,.
$$
%Taking into account the definition \eqref{sec:Ad.4-1} we set 
Setting now
$\wt{\omega}=\omega(\wt{\alpha})$ with the function $\omega$ defined in \eqref{sec:Ad.4-1},
the index $\wt{\alpha}$ defined in \eqref{sec:Up.0-1},
$\wt{j}_0=[\wt{\omega} \ve_\zs{n}]$,
$\wt{j}_1=[\wt{\omega}/\ve_\zs{n}]$
and
$$
\varsigma_\zs{n}(S)=\frac{1}{n}\,
\sum^n_\zs{l=1}\sigma^2_\zs{l}(S)\,,
$$
we rewrite \eqref{sec:Up.1-1} as follows
\begin{align}\label{sec:Up.1-2}
\E_\zs{S,p}\,\|\wt{S}-S\|_\zs{n}^2&=\sum_{j=\wt{j}_0+1}^{\wt{j}_1-1}
(1\,-\,\wt{\lambda}_\zs{j})^2\theta^2_\zs{j,n}\\ \nonumber
&+\frac{\varsigma_\zs{n}(S)}{n}\,
\sum_{j=1}^n\,\wt{\lambda}_\zs{j}^2+\wt{\Delta}_\zs{1,n}+\wt{\Delta}_\zs{2,n}
\end{align}
with 
$$
\wt{\Delta}_\zs{1,n}\,=\,\sum_{j=\wt{j}_1}^{n}\,\theta^2_\zs{j,n}
\quad
\mbox{and}
\quad
\wt{\Delta}_\zs{2,n}\,=\frac{1}{n}\,\sum_{j=1}^n\,\wt{\lambda}^2_\zs{j}
\left(\varsigma_\zs{j,n}(S)-\varsigma_\zs{n}(S)\right)\,.
$$
Note that we have decomposed the first term on the right-hand of \eqref{sec:Up.1-1} into the sum
$$
\sum_{j=\wt{j}_0+1}^{\wt{j}_1-1}
(1\,-\,\wt{\lambda}_\zs{j})^2\theta^2_\zs{j,n}\,+\,\wt{\Delta}_\zs{1,n}\,.
$$
This decomposition allows us to show that $\wt{\Delta}_\zs{1,n}$ is negligible and further to
approximate the first term by a similar term in which the coefficients $\theta_\zs{j,n}$
will be replaced by the Fourier coefficients $\theta_\zs{j}$ of the function $S$.

Taking into account the definition of $\omega$ in \eqref{sec:Ad.4-1},
we can bound $\wt{\omega}$ as
$$
\wt{\omega} \ge (A_\zs{k})^{\frac{1}{2k+1}}\,
(n \ve_\zs{n})^{\frac{1}{2k+1}}\,.
$$
Therefore, by
Lemma~\ref{Le.A.1} we obtain 
$$
\lim_{n\to\infty}\sup_{S\in W_r^k}\,n^{\frac{2k}{2k+1}}\,\wt{\Delta}_\zs{1,n}=0\,.
$$
Let us consider now the next term $\wt{\Delta}_\zs{2,n}$. We have
\begin{align*}
|\wt{\Delta}_\zs{2,n}|
=\frac{1}{n^2}\,\left|\sum^n_{i=1}\,\sigma^2_\zs{i}(S)\,
\sum_{j=1}^n\,\wt{\lambda}^2_\zs{j}\,
(\phi^2_\zs{j}(x_\zs{i})-1) \right|
\le \frac{\sigma_\zs{*}}{n}
\sup_\zs{0\le x\le 1}
\left|\sum_{j=1}^n\,\wt{\lambda}_\zs{j}^2\,
(\phi^2_\zs{j}(x)-1)\right|\,.
\end{align*}
%where $\ov{\phi}_\zs{j}(x)=\phi^2_\zs{j}(x)-1$. 
Now by Lemma~\ref{Le.A.2} and the
definition \eqref{sec:Ad.4},
we obtain directly the same property for $\wt{\Delta}_\zs{2,n}$, i.e.
$$
\lim_{n\to\infty}\sup_{S\in W_r^k}\,n^{\frac{2k}{2k+1}}\,|\wt{\Delta}_\zs{2,n}|=0\,.
$$
Setting 
$$
\wh{\gamma}_\zs{k,n}(S)=n^{\frac{2k}{2k+1}}
\sum_{j=\wt{j}_0}^{\wt{j}_1-1}(1-\wt{\lambda}_\zs{j})^2\theta^2_\zs{j}
+
\varsigma_\zs{n}(S) n^{-\frac{1}{2k+1}}
\sum_{j=1}^n\,\wt{\lambda}_\zs{j}^2
$$
and applying the well-known inequality
$$
(a+b)^2\le (1+\delta)a^2+(1+1/\delta)b^2
$$
to the first term on the right-hand side of the inequality \eqref{sec:Up.1-2}, we obtain that, for any  $\delta>0$
and for any $p\in\cP_\zs{n}$,
\begin{align}\nonumber
\E_\zs{S,p}\,\|\wt{S}-S\|_\zs{n}^2 &\le(1+\delta)
\,\wh{\gamma}_\zs{k,n}(S)\,n^{-2k/(2k+1)}\\[2mm] \label{sec:Up.1-3}
&+\wt{\Delta}_\zs{1,n}+
\wt{\Delta}_\zs{2,n}+(1+1/\delta)\,
\wt{\Delta}_\zs{3,n}\,,
\end{align}
where 
$$
\wt{\Delta}_\zs{3,n}=\sum_{j=\wt{j}_0+1}^{\wt{j}_1-1}(\theta_\zs{j,n}\,-\,\theta_\zs{j})^2\,.
$$
Taking into account 
that $k\ge 1$ and that
$$
\wt{j}_\zs{1}\le 
\ov{\omega}\,\ve^{-1}_\zs{n}
+(A_\zs{k})^{\frac{1}{2k+1}}\,
n^{\frac{1}{2k+1}} 
(\ve_\zs{n})^{-(2k+2)/(2k+1)}\,,
$$
we can show through Lemma~\ref{Le.A.3} 
that
$$
\lim_{n\to\infty}\sup_{S\in W_r^k}\,n^{\frac{2k}{2k+1}}\,\wt{\Delta}_\zs{3,n}=0\,.
$$
Therefore, the inequality  \eqref{sec:Up.1-3}
 yields
$$
\limsup_{n\to\infty} n^{\frac{2k}{2k+1}}
\sup_{S\in W^k_r}
\cR_\zs{n}(\wt{S},S)/\gamma_k(S)\le 
\limsup_{n\to\infty}
\sup_{S\in W^k_r}
\wh{\gamma}_\zs{k,n}(S)/\gamma_k(S)
$$
and to prove \eqref{sec:Up.1} it suffices to show that
\begin{equation}\label{sec:Up.3}
\limsup_{n\to\infty}
\sup_{S\in W^k_r}
\wh{\gamma}_\zs{k,n}(S)/\gamma_k(S)
\le 1\,.
\end{equation}
First, it should be noted 
that the definition \eqref{sec:Up.0}
and the
inequalities \eqref{sec:Co.6}-\eqref{sec:Co.6-1} imply directly
$$
\lim_\zs{n\to\infty}\sup_{S\in W^k_r}
\left|
\wt{t}_\zs{n}/\ov{r}(S)
-1
\right|=0\,.
$$
Moreover,
by the definition of $(\wt{\lambda}_\zs{j})_\zs{1\le j\le n}$ in 
\eqref{sec:Up.0-2},
for sufficiently large $n$ for which $\wt{t}_\zs{n}\ge \ov{r}(S)$,
we find
\begin{align*}
\sup_\zs{j\ge 1}
\,n^{\frac{2k}{2k+1}}
\,
\frac{(1-\wt{\lambda}_\zs{j})^2}{(\pi j)^{2k}}
&= \pi^{-2k}
(A_\zs{k}\wt{t}_\zs{n})^{-2k/(2k+1)}
\le \pi^{-2k}
(A_\zs{k}\ov{r}(S))^{-2k/(2k+1)}\,.
\end{align*}
 Therefore, by the definition of the coefficients $(a_\zs{j})_\zs{j \ge 1}$
in \eqref{sec:Co.4}, one has 
$$
\limsup_{n\to\infty}
n^{\frac{2k}{2k+1}}
\sup_{S\in W^k_r}
\sup_{j\ge \wt{j}_0}\,\pi^{2k}
(A_\zs{k}\ov{r}(S))^{2k/(2k+1)}
(1-\wt{\lambda}_\zs{j})^2/a_\zs{j}
\le 1\,.
$$
Furthermore, in view of the definition \eqref{sec:Ad.4} we calculate directly
$$
\lim_{n\to\infty}
\sup_{S\in W^k_r}
\left|
\,
n^{-\frac{1}{2k+1}}\,\sum^n_\zs{j=1}\wt{\lambda}^2_\zs{j}
-
(A_\zs{k}\ov{r}(S))^{\frac{1}{2k+1}}
\int^1_\zs{0}(1-z^k)^2\d z
\right|=0\,.
$$
Now, the definition of $W^k_\zs{r}$ in \eqref{sec:Co.3}
and the condition \eqref{sec:Co.5} imply the inequality \eqref{sec:Up.3}. Hence Theorem~\ref{Th.sec:Up.1}.
\endproof

\subsection{ Proof of Theorem~\ref{Th.sec:Tr.1}}\label{subsec:A.2}

For any $z=(z_\zs{1},\ldots,z_\zs{l})'\in\bbr^n$, we set 
$$
\varrho_\zs{i}(v,z)=\frac{1}{f(v,z)\Phi(z)}
\,\frac{\partial}{\partial z_\zs{i}}\left(f(v,z)\Phi(z)\right)\,.
$$
Note that, due to the condition \eqref{sec:Co.6}, the density \eqref{subsec:Tr.2} is bounded, i.e.
$$
f(v,z)\le (2\pi g_\zs{*})^{-n/2}\,.
$$
So through \eqref{subsec:Tr.1-1} we obtain that
$$
\lim_\zs{|z_\zs{i}|\to\infty}
\tau(z)\,f(v,z)\varphi_\zs{i}(z_\zs{i})=0\,.
$$
Therefore, integrating by parts yields
\begin{align*}
\wt{\E}(\wh{\tau}_\zs{n}-\tau(\vartheta))\varrho_\zs{i}&=
\int_\zs{\bbr^{n+l}}
(\wh{\tau}_\zs{n}(v)-\tau(z))
\frac{\partial}{\partial z_\zs{i}}
\left(f(v,z)\Phi(z)\right)\d z\,\d v\\
&=
\int_\zs{\bbr^l}\left(\frac{\partial}{\partial\, z_\zs{i}}\tau(z)\right)
\Phi(z)\left(\int_\zs{\bbr^{n}} f(v,z)\d v \right)\d z=\ov{\tau}_\zs{i}\,.
\end{align*}
Now the Bounyakovskii-Cauchy-Schwarz inequality gives the following lower bound
$$
\wt{\E}(\wh{\tau}_\zs{n}-\tau(\vartheta))^2\ge
\ov{\tau}_\zs{i}^2/\wt{\E}\varrho_\zs{i}^2\,.
$$
To estimate the denominator in the last ratio, note that
$$
\varrho_\zs{i}(v,z)
=\wt{f}_\zs{i}(v,z)
+
\frac{\dot{\varphi}_\zs{i}(z_\zs{i})}{\varphi_\zs{i}(z_\zs{i})}
\quad\mbox{with}\quad
\wt{f}_\zs{i}(v,z)=
(\partial/\partial\,z_\zs{i})\ln f(v,z)\,.
$$
From \eqref{subsec:Tr.1} it follows that
$$
\wt{f}_\zs{i}(v,z)=
\sum^{n}_\zs{j=1}(\xi^2_\zs{j}-1)\,
\frac{1}{2\sigma^2_\zs{j}(z)}
\frac{\partial}{\partial\,z_\zs{i}}\,\sigma^2_\zs{j}(z)
+
\sum^{n}_\zs{j=1}\,\xi_\zs{j}\,
\frac{S^{\prime}_\zs{i}(x_\zs{j})}{\sigma_\zs{j}(z)}\,.
$$
Moreover, the conditions $\H_\zs{2})$ and \eqref{subsec:Tr.2-1} imply
\begin{align*}
(\partial/\partial\,z_\zs{i})\,\sigma^2_\zs{j}(z)\,=\,
(\partial/\partial\,z_\zs{i})\,g^2(x_\zs{j},S_\zs{z})\,=\,
\wt{\L}_\zs{i}(x_\zs{j},
z)
\end{align*}
from which it follows
$$
\wt{\E}\,\left(
\wt{f}_\zs{i}(Y,\vartheta)
\right)^2\,
=F_\zs{i}+B_\zs{i}\,. 
$$
This implies inequality \eqref{subsec:Tr.3}. Hence Theorem~\ref{Th.sec:Tr.1}.
\endproof
\medskip
\subsection{ Proof of Theorem~\ref{Th.sec:Ma.3}}\label{subsec:Th.M.3}

In this section we prove Theorem~\ref{Th.sec:Ma.3}.
 Lemma~\ref{Le.sec:Up.1} implies that to prove the lower bounds \eqref{sec:Ma.6} and \eqref{sec:Ma.6-1}, it suffices
to show 
\begin{equation}\label{sec:Lo.1}
\liminf_{n\to\infty}\,\inf_{\wh{S}_\zs{n}}\,n^{\frac{2k}{2k+1}}\,
\cR_0(\wh{S}_\zs{n})\,\ge\,1\,,
\end{equation}
where
$$
\cR_0(\wh{S}_\zs{n})\,=\,
\sup_\zs{S\in W_r^k}\,\E_\zs{S,q}\,\|\wh{S}_\zs{n}-S\|^2/\gamma_k(S)\,.
$$

 For any estimator $\wh{S}_\zs{n}$, we denote by $\wh{S}^0_\zs{n}$ its projection onto $W_r^k$, i.e.\\
 $\wh{S}^0_\zs{n}=\hbox{\rm Pr}_\zs{W_r^k}(\wh{S}_\zs{n})$.
Since $W^k_\zs{r}$ is a convex set, we get 
$$
\|\wh{S}_\zs{n}-S\|^2\ge\|\wh{S}^0_\zs{n}-S\|^2\,.
$$
Now  we introduce the following set
\begin{equation}\label{Fa.11}
\Xi_\zs{n}=\{\max_\zs{1\le m\le M_\zs{n}}\,\max_\zs{1\le j\le N}\,\zeta^2_\zs{m,j}\le d_\zs{n}\}\,,
\end{equation}
where  $(\zeta_\zs{m,j})$ are i.i.d. $\cN(0,1)$ random variables 
from \eqref{Fa.5}
and the sequence $(d_\zs{n})_\zs{n\ge 1}$ is given in the condition $\A_\zs{2})$.
It is clear that  the condition $\A_\zs{1})$  implies
\begin{equation}\label{Fa.12}
\limsup_\zs{n\to\infty}\,n^\nu\,
\P\left(\Xi^c_\zs{n}\right)=0\,.
\end{equation}

Therefore, we can write that
$$
\cR_0(\wh{S}_\zs{n})
\ge\int_{\{z:S_\zs{z,n}\in W^k_\zs{r}\}\cap\Xi_\zs{n}}
\,\frac{\E_\zs{S_\zs{z,n},q}\|\wh{S}^0_\zs{n}-S_\zs{z,n}\|^2}{\gamma_k(S_\zs{z,n})}\,\mu_\zs{\vartheta}(\d z)\,.
$$
Here the kernel function family $(S_\zs{z,n})$ is given in \eqref{Fa.4}
in which $N_\zs{n}=[\ln^4 n]+1$ and
 the parameter $h$ is defined in \eqref{sec:Lo.17} and \eqref{sec:Lo.20}; the measure
 $\mu_\zs{\vartheta}$ is defined in \eqref{sec:Lo.3}.
Moreover, note that on the set $\Xi$ the random function $S_\zs{\vartheta,n}$
is uniformly bounded, i.e.
\begin{equation}\label{Fa.18}
\|S_\zs{\vartheta,n}\|_\zs{\infty}=
\sup_\zs{0\le x\le 1}\,|S_\zs{\vartheta,n}(x)|\,
\le \sqrt{d_\zs{n}}\,t^*_\zs{n}\,,
\end{equation}
where the coefficient $t^*_\zs{n}$ is defined in \eqref{Fa.6}.
 
Thus, we estimate the risk $\cR_0(\wh{S}_\zs{n})$ from below as
$$
\cR_0(\wh{S}_\zs{n})
\ge
\frac{1}{\gamma^*_\zs{n}}\,
\int_{\{z:S_\zs{z,n}\in W^k_\zs{r}\}\cap\Xi_\zs{n}}\,
\E_\zs{S_\zs{z,n},q}\|\wh{S}^0_\zs{n}-S_\zs{z,n}\|^2\,\mu_\zs{\vartheta}(\d z)
$$
with 
\begin{equation}\label{sec:Lo.2}
\gamma^*_\zs{n}=\sup_\zs{\|S\|_\zs{\infty}\le \sqrt{d_\zs{n}}t^*_\zs{n}}\,\gamma_k(S)\,.
\end{equation}

By making use of the  Bayes risk \eqref{sec:Lo.3} with the prior distribution given by 
formulas
 \eqref{Fa.5}, \eqref{sec:Lo.15-1}, \eqref{sec:Lo.17}
and \eqref{sec:Lo.20}, for any fixed parameter $0<\ve<1$,
 we  rewrite  the lower bound for $\cR_0(\wh{S}_\zs{n})$ as 
\begin{equation}\label{sec:Lo.4}
\cR_0(\wh{S}_\zs{n})
\ge
\cE_\zs{\ve,n}(\wh{S}^0_\zs{n})/\gamma^*_\zs{n}-2\,
\Omega_\zs{n}/\gamma^*_\zs{n}
\end{equation}
with
$$
\Omega_\zs{n}=\E
(\Chi_\zs{\{S_\zs{\vartheta,n}\notin W^k_\zs{r}\}}\,+\,
\Chi_\zs{\Xi^c_\zs{n}})
(r+\|S_\zs{\vartheta,n}\|^2)\,.
$$

In Section~\ref{subsec:Br} we proved that the parameters in chosen prior distribution 
 satisfy the conditions $\A_\zs{1})$--$\A_\zs{4})$.
Therefore, Propositions~\ref{sec:Pr.Fa.2}--\ref{sec:Pr.Fa.3} and the limit \eqref{Fa.12}
imply that, for any $\nu>0$,
$$
\lim_\zs{n\to\infty}\,n^{\nu}\,
\Omega_\zs{n}=0\,.
$$
Moreover, by the condition $\H_\zs{4})$ the sequence $\gamma^*_\zs{n}$ goes to $\gamma_k(S_\zs{0})$
as $n\to\infty$. Therefore, from this, \eqref{sec:Lo.21} and \eqref{sec:Lo.4} we get, for any $0<\ve<1$,
$$
\liminf_{n\to\infty}\,\inf_{\wh{S}_\zs{n}}\,n^{\frac{2k}{2k+1}}\,
\cR_0(\wh{S}_\zs{n})\,\ge\,\frac{(1+\ve^\prime)(1-\ve)^{\frac{1}{2k+1}}}{(1+\ve)^{\frac{1}{2k+1}}}\,.
$$
where $\ve^\prime$ is defined in \eqref{sec:Lo.19-1}.
Limiting here $\ve\to 0$ implies the inequality \eqref{sec:Lo.1}. 
Hence Theorem~\ref{Th.sec:Ma.3}.
\endproof

\vspace{10mm}
{\bf Acknowledgements}

We are grateful for the comments and stimulating questions by the associate Editor and two anonymous Referees which have led to considerable improvements.

%\newpage

\begin{flushright}
\begin{tabular}{lcl}
   L.Galtchouk                       &$\quad$& S. Pergamenshchikov              \\
 Department of Mathematics           &$\quad$& Laboratoire de Math\'ematiques Raphael Salem,\\      
 Strasbourg University               &$\quad$& Avenue de l'Universit\'e, BP. 12,            \\
 7, rue Rene Descartes               &$\quad$&  Universit\'e de Rouen,                  \\
 67084, Strasbourg, France           &$\quad$&  F76801, Saint Etienne du Rouvray, Cedex France.\\
 e-mail: galtchou@math.u-strasbg.fr  &$\quad$& Serge.Pergamenchtchikov@univ-rouen.fr         \\
\end{tabular}
\end{flushright} 
%\end{description}

 \end{document}